\documentclass[a4,11pt]{article}
\usepackage{amssymb}
\usepackage{amsmath}
\usepackage{bbm}
\usepackage{empheq}
\usepackage[framemethod=tikz]{mdframed}
\usepackage{lipsum}
\usepackage{tcolorbox}
\usepackage{fourier-orns}
\usepackage{enumitem}
\RequirePackage[amsmath,thmmarks,hyperref,framed]{ntheorem}

\usepackage[applemac]{inputenc}

\date{}

\newtheorem{theorem}{Theorem}
\newtheorem{prop}[theorem]{Proposition}
\newtheorem{cor}[theorem]{Corollary}

\newtheorem{lemma}[theorem]{Lemma}

\numberwithin{equation}{section}

\newcommand{\apropto}{\mathrel{\vcenter{
  \offinterlineskip\halign{\hfil$##$\cr
    \propto\cr\noalign{\kern2pt}\sim\cr\noalign{\kern-2pt}}}}}

\usepackage[english]{babel}

\theorembodyfont{\upshape}

\title{Precise asymptotics for the density and the upper tail of exponential functionals of subordinators}

\author{B\'en\'edicte Haas\thanks{Universit\'e Sorbonne Paris Nord, LAGA, CNRS (UMR  7539) 93430 Villetaneuse, France \newline \hspace*{0.5cm} E-mail: haas@math.univ-paris13.fr}}

\begin{document}

\maketitle

\begin{abstract}
We provide exact large-time equivalents of the density and upper tail distributions of the exponential functional of a subordinator in terms of its Laplace exponents. This improves previous results on the logarithmic asymptotic behaviour of the upper tail. Since the results in their general form remain rather abstract, we study some fairly common situations in details to see how the behavior near 0 of the L\'evy measure of the subordinator influences the asymptotics of the density and upper tail distributions. Several examples of applications are then discussed. 
\end{abstract}

\section{Introduction}

The exponential functional of a subordinator $\xi$ is the random variable
$$
I:=\int_0^{\infty} \exp(-\xi_r) \mathrm dr.
$$
To simplify we always assume in this paper that $\xi$ starts from 0. Such random variables have been thoroughly investigated, since they play a key role in several facets of modern probabilities since they correspond to the extinction times of non-increasing self-similar Markov processes by Lamperti's transformation \cite{Lamp72}. As such, they are involved in the description of various processes ranging from the analysis of algorithms to coagulation or fragmentation processes \cite{BertoinSSF,DGR02,GRZ04,IM08,HM11}. They belong to the more general class of exponential functionals of L\'evy processes drifting to $\infty$, which offers a wider range of applications. We refer to the survey of Bertoin and Yor \cite{BYSurvey} for a general overview on that topic.

\bigskip

A well-known result of Carmona, Petit and Yor \cite[Proposition 2.1]{CPY97}, generalized to any subordinator by Pardo, Rivero and Schaik \cite[Theorem 2.3]{PRvS13}, is that the law of $I$ has a density with respect to the Lebesgue measure on $[0,\infty)$, which satisfies an integral equation recalled in (\ref{eq:k}). Like most articles on that topic, we denote this density by $k$. This density, or the distribution of $I$, are known explicitly only in a few cases. For instance, \cite{BYSurvey} identifies subordinators that have an exponential functional $I$ distributed as the power of an exponential distribution or as a generalized Mittag-Leffler distribution. Fine information is also known  in some other specific cases, for instance when $\xi$ is a standard Poisson process \cite{BBY04,DGR02} or more generally a compound Poisson process with jumps having a fractional generating function \cite{GRZ04}, a series representation of $k$ is obtained explicitly using $q$-calculus techniques.

\bigskip

In this paper we obtain an equivalent of the density $k(t)$ and the (upper) tail $\mathbb P(I>t)$ as $t \rightarrow \infty$ in a general framework. We also obtain the first order of convergence and give the tools to obtain further orders. It is well known that the tail of $I$ is exponential or even lighter than exponential. More precisely a first order equivalent of the asymptotic behavior of the \emph{logarithm} of the tail has been obtained by Rivero \cite{Rivero03} when the L\'evy measure of the subordinator is infinite whereas  Maulik and Zwart \cite{MZ06} studied the asymptotic behavior of the tail in the compound Poisson case. Their results will be recalled below, once we have introduced the necessary notations.  For comparison, the tail of the exponential functional of a non-monotonic L\'evy process shows a drastically different asymptotic behavior since it is generally heavy tailed or even extremely heavy tailed (it may decrease as the power of the logarithm): see again Maulik and Zwart \cite{MZ06} and Rivero \cite{R05,R12} for results in a general framework, and Kuznetsov and Pardo \cite{K12,KP13}, Patie \cite{Patie12} and the references therein for series representation of $k$ for some specific subclasses of non-monotonic L\'evy processes. To complete this picture we mention that some results on the small-time behavior of the density and the tail are also available in the literature, see e.g. \cite{PRvS13, PS18} and the references therein.

\bigskip

\textbf{Examples and applications.} The expressions of the equivalents of the density and tail of $I$ remain quite abstract in their general form and it is interesting to see on examples what is the precise influence of the behavior of the L\'evy measure near 0. Even when this L\'evy measure is finite (compound Poisson cases),  one observes a wide range of behaviors (see Case 1 in Section \ref{sec:MR}). We detail in this paper several examples, typically when the lower and upper tails of the L\'evy measure have asymptotic power expansions at 0, but not only, and then discuss several applications to non-negative non-decreasing self-similar Markov processes, and fragmentation and coalescence processes.

\bigskip

\textbf{Organization of the paper.} Our main contribution on the asymptotic behavior of the tail and density of $I$ is given in Section \ref{sec:MRA}, where we also present the examples and applications. Section \ref{sec:phi,psi} is devoted to some elementary properties of the functions $\phi$ and $\psi$ which will be needed for the proof of our main result, Theorem \ref{thm:main}, undertaken in Section \ref{sec:thmmain}. Next, we focus on the applications: Section \ref{sec:special} details the cases where the Laplace exponent $\phi$ has an asymptotic power expansion, while Section \ref{sec:ssMp} provides a proof of Corollary \ref{cor:ssMp} on the asymptotics of moments of non-increasing, non-negative self-similar Markov processes. 

\bigskip

\textbf{Note added.} While completing the writing of this manuscript, we became aware that the equivalents of the density and upper tail of $I$ are also studied by Minchev and Savov in a paper posted recently on the arXiv \cite{MS21}. Our approaches are different and each seems to have its advantages. Minchev and Savov apply the saddle point method to the Mellin transform of the exponential functional, while we exploit the proximity of an integral equation satisfied by the derivative of the logarithm of the tail of $I$ with an integral equation satisfied by a function related to the Laplace exponent of the subordinator. This allows us to get immediately the first order of convergence and a method to give any further order by iterating our calculations. On the other hand, Minchev and Savov get an expression of a constant of multiplicity (noted $c_I$ in the rest of this article) in terms of the Laplace exponent of the subordinator, which we do not have. They are also interested in the equivalents of the derivatives of the density. Here we also develop several examples and applications.

\section{Main results and applications}
\label{sec:MRA}

\subsection{Background and notation} 

Throughout we assume that the subordinator $\xi$ is drift-free, unkilled and start from 0. We call $\pi$ its L\'evy measure (which is assumed to be non-zero) and $\phi$ its Laplace exponent, so that $\mathbb E[\exp(-x\xi_t)]=\exp\left(-t \phi(x)\right), \forall t,x \geq 0$ with 
\begin{equation}
\label{def:phi}
\phi(x):=\int_0^{\infty} \left(1-\exp(-xv)\right)\pi(\mathrm dv), \qquad  x \geq 0.
\end{equation}
When the L\'evy measure $\pi$ is infinite and more precisely when the function $\phi$ is regularly varying at $\infty$ with some index $\gamma \in (0,1)$,  Rivero \cite[Proposition 2]{Rivero03} proved that 
\begin{equation*}
\label{eq:Rivero}
\ln \mathbb P(I>t)  \underset{t \rightarrow \infty} \sim  -(1-\gamma) \psi(t),
\end{equation*}
where the function $\psi$ is the inverse of the continuous and strictly increasing function $x \mapsto x/\phi(x)$:
\begin{equation}
\label{def:psi}
\frac{\psi(x)}{\phi(\psi(x))}=x.
\end{equation}
Note that $\psi$ is well-defined on $(x_{\psi},\infty)$ where $$x_{\psi}:=\lim_{x \downarrow 0} x/\phi(x)=\left(\int_0^{\infty} v \pi(\mathrm dv)\right)^{-1}$$ is finite (possibly null), and that $\phi,\psi$ are infinitely differentiable on $(0,\infty)$ and $(x_{\psi},\infty)$ respectively. In \cite[Proposition A.1]{HR12} an extension of Rivero's result to the case where 
\begin{equation}
\label{hyp:main}
\tag{$\mathbf H$}
\limsup_{x \rightarrow \infty} \frac{x\phi'(x)}{\phi(x)} <1
\end{equation}
is given; then 
\begin{equation}
\label{extR}
\ln \mathbb P(I>t)  \underset{t \rightarrow \infty}\sim -\int_{x_{\psi}}^{t} \frac{\psi(r)}{r} \mathrm dr. 
\end{equation}

\smallskip

In this paper we will always work under this assumption (\ref{hyp:main}). We highlight that it is weaker than assuming that $\phi$ is regularly varying at $\infty$, hence in particular that $\pi$ is finite. We note also that without any hypothesis on the model, $\phi'(x)\leq \phi(x)/x$ for all $x \geq 0$ since $\phi$ is concave and thus the above limsup is always smaller or equal to 1. The assumption  (\ref{hyp:main}) is therefore not very restrictive and in fact verified in all applications we are interested in.  

\medskip

When $\pi$ is finite Maulik and Zwart \cite[Sections 5 and 6]{MZ06} obtained more precise results, which can be summarized as follows:

\smallskip
\begin{itemize}[topsep=0pt,itemsep=5pt,parsep=0pt,partopsep=0pt,leftmargin=15pt]
\item[$\bullet$] when $\int_0^{\infty} u^{-1} \pi(\mathrm du)<\infty$,  $$\mathbb P(I>t)  \underset{t \rightarrow \infty}\sim \mathrm c_I \exp(-|\pi|t)$$ where $|\pi|=\pi(0,\infty)$, $\mathrm c_I=\mathbb E[e^{|\pi|e^{-X}I}]$ is finite and positive, and $X$ has the distribution of a jump of the subordinator ($\pi/|\pi|$) and is independent of $I$
\item[$\bullet$] when $\pi(0,u) =bu+o(u^{1+\delta})$ as $u \downarrow 0$ for some $b,\delta>0$,  $$\mathbb P(I>t)  \underset{t \rightarrow \infty}\sim \mathrm c_I ~ t^{\frac{b}{|\pi|}}\exp(-|\pi|t)$$ where $\mathrm c_I=|\pi|^be^{b \gamma_{\mathrm{Euler}}}\left(\prod_{k=1}^{\infty}(1-\int_0^{\infty}e^{-uk}\pi(\mathrm du)/|\pi|)e^{b/k}\right)^{-1}$ is finite and positive, and $\gamma_{\mathrm{Euler}}$ is Euler's constant.
\end{itemize}

\bigskip

\textbf{Notation.} Throughout the paper, we will use the standard big-O and little-o Landau asymptotic notation, as well as  the notation $f\lesssim g$ on a domain $D \subset \mathbb R^d$, meaning that if $f,g$ are two positive functions defined on $D$, there exists a positive constant $\mathrm c$ such that $f(x) \leq \mathrm c g(x)$, $\forall x \in D$. Last, for two functions $f,g$ defined on an interval $[a,\infty)$, we will write $f(x) ~ \underset{x \rightarrow \infty}{\sim}~ g(x)$ when $f(x)/g(x) \underset{x \rightarrow \infty}{\rightarrow} 1$ and more generally $f(x) ~ \underset{x \rightarrow \infty}{\propto}~ g(x)$ when $f(x)/g(x) \underset{x \rightarrow \infty}{\rightarrow} \mathrm c$ for some $\mathrm c>0$.

\subsection{Main results} 
\label{sec:MR}

Our main contribution is the following theorem and its corollary which improve Rivero's result by giving an exact (non-logarithmic) asymptotic behavior of the tail of $I$, as well as its density $k$, and completes Maulik and Zwart' results in the compound Poisson case.

\begin{theorem}
\label{thm:main}
Assume \emph{(\ref{hyp:main})}. Then there exists a constant $\mathrm c_I \in (0,\infty)$ such that as $t \rightarrow \infty$
$$
\mathbb P(I>t) ~ =  ~ \mathrm c_I ~  \frac{t (\psi'(t))^{1/2}}{\psi(t)} ~ \exp\left(- \int_{x_{\psi}+1}^t \frac{\psi(r)}{r} \mathrm dr\right)\left( 1+O\left(\frac{1}{\psi(t)}\right)\right)
$$
and
$$
k(t) ~ =  ~  \mathrm c_I ~ \left(\psi'(t)\right)^{1/2} ~ \exp\left(- \int_{x_{\psi}+1}^t \frac{\psi(r)}{r} \mathrm dr\right) \left( 1+O\left(\frac{1}{\psi(t)}\right)\right).
$$
\end{theorem}

\bigskip

We obtain in fact an asymptotic expansion of the quotient $k(t)/\mathbb P(I>t)$ up to a $O(1/(\psi(t))^2)$ term in Proposition \ref{prop:equiv}, which is the key technical result of this paper and implies directly the above theorem. This is proved by using iteratively the integral equation (\ref{eq:k}) satisfied by $k$, which permits us to evaluate how close the quotients $k(t)/\mathbb P(I>t)$ and $\psi(t)/t$ are as $t \rightarrow \infty$. By pushing forward the iterative argument we could in fact improve the asymptotic expansion of Proposition \ref{prop:equiv} to any order. \emph{This would allow us to develop the parentheses $(1+O(1/\psi(t)))$ of Theorem \ref{thm:main} to any order}. However this does not seem to give easy-to-state expressions and we will be content with Proposition \ref{prop:equiv} and Theorem \ref{thm:main} here, and leave the additional calculations to the interested reader if they need further orders.  Compared to previous works, we note that several results recalled at the beginning of this introduction were obtained by inverting the Mellin or Laplace transforms of the distribution of $I$, a technique that we do not use here. In spirit, we are closer to the approach of Maulik and Zwart \cite{MZ06} for finite measures $\pi$ which relies on the fact that the distribution of $I$ then satisfies the random affine equation  ~$I=E+\exp(-X) \tilde I$~ where $E,X,\tilde I$ are independent random variables, $E$ having an exponential distribution, $X$ being distributed according to $\pi/|\pi|$ and $\tilde I $ as $I$. However for infinite L\'evy measure $\pi$ there is no such explicit equation satisfied by $I$ in general.

\bigskip

In the integral of Theorem \ref{thm:main} one can of course replace the bound $x_{\psi}+1$ by any $x \in (x_{\psi},\infty)$, by multiplying the constant $\mathrm c_I$ appropriately. 
Our approach does not allow us, a priori, to get an explicit expression of $\mathrm c_I$. See the recent paper \cite{MS21} for an expression of this constant involving an integral of a function of $\phi$ and its first and second derivatives.

\bigskip

Note that if $I(\pi)$ denotes the exponential functional of a subordinator with L\'evy measure $\pi$ then $I(c\pi)\overset{\mathrm{(d)}}=c^{-1}I(\pi)$ for any $c>0$. So in the examples we will treat, we will always try to take the multiplicative version of the L\'evy measure that gives the easiest to read asymptotic expressions of the tail and density of $I$. 

\bigskip

There is one useful and immediate corollary of the above theorem. Indeed, when $\phi$ is regularly varying with some index $\gamma \in [0,1)$ one easily checks that $x \psi'(x)/\psi(x) \rightarrow (1-\gamma)^{-1}$ as $x \rightarrow \infty$ (see Lemma \ref{cororegvar}) and the above theorem then implies:

\bigskip

\begin{cor}[Regular variation]
\label{cor:RV}
Assume that $\phi$ is regularly varying at $\infty$ with some index $\gamma \in [0,1)$. Then 
$$
\mathbb P(I>t) ~ \underset{t \rightarrow \infty}{\sim}~ \frac{\mathrm c_I}{(1-\gamma)^{1/2}} ~ \left(\frac{t}{\psi(t)}\right)^{1/2} \exp\left(- \int_{x_{\psi}+1}^t \frac{\psi(r)}{r} \mathrm dr\right)
$$
and
$$
k(t) ~  ~ \underset{t \rightarrow \infty}{\sim}  ~  \frac{\mathrm c_I}{(1-\gamma)^{1/2}} ~  \left(\frac{\psi(t)}{t}\right)^{1/2} \exp\left(- \int_{x_{\psi}+1}^t \frac{\psi(r)}{r} \mathrm dr\right),
$$
where $\psi(t)/t$ is regularly varying at $\infty$ with index $\gamma/(1-\gamma)$.
\end{cor}

\bigskip

In practice, to extract more explicit asymptotic expressions from those of Theorem \ref{thm:main} and Corollary \ref{cor:RV}, the general strategy is of course to find a sufficiently precise asymptotic expansion for $\phi$ and then deduce one for $\psi$. To illustrate, we finish this section by detailing cases where the lower or upper tails of $\pi$ have asymptotic power expansions. We refer to Section \ref{sec:special} for some details on calculations. 

\bigskip

\textbf{Case 1.} Assume that $\pi$ is \emph{finite}, say $\pi$ is a probability. Then, $x^{-1}\psi(x)=\phi(\psi(x))=1-\int_0^{\infty}\exp(-\psi(x)u) \pi(\mathrm du)$ which shows both that $x^{-1}\psi(x) \rightarrow 1$ and then $x^{-1}\psi(x)=$ \linebreak $1+O \left(\int_0^{\infty} \exp(-xu/2)\ \pi(\mathrm du)\right) $ as $x \rightarrow \infty$. We need more information on $\pi$ to understand the behavior of the Big-O. Some examples: 

\smallskip

\begin{itemize}[topsep=0pt,itemsep=5pt,parsep=0pt,partopsep=0pt,leftmargin=15pt]
\item[$\bullet$] If $\int_0^{\infty} u^{-1}\pi(\mathrm du)<\infty$, then the above remark shows that $\int_{x_{\psi}+1}^t r^{-1}\psi(r) \mathrm dr-t$ converges as $t \rightarrow \infty$ and we get that
$$
\mathbb P(I>t)  \underset{t \rightarrow \infty }{\sim}  k(t) \underset{t \rightarrow \infty }{\propto}  \exp(-t).
$$
\item[$\bullet$] Now, assume that $\pi(0,u)=\sum_{i=1}^pc_iu^{\gamma_i}+O(u^{1+\varepsilon})$ as $u\rightarrow 0$ where $0<\gamma_1<\ldots<\gamma_{p-1}<\gamma_p=1$, $\varepsilon>0$ and the $c_i$s may be any real numbers. Then, if moreover $\gamma_1>1/2$,
\begin{eqnarray*}
\mathbb P(I>t)  \underset{t \rightarrow \infty }{\sim}  k(t)
 \underset{t \rightarrow \infty }{\propto}   t^{c_p} \exp\left(- t+\sum_{i=1}^{p-1} \frac{c_i\Gamma(1+\gamma_i)}{1-\gamma_i} t^{1-\gamma_i}\right).
\end{eqnarray*}
If some of the $\gamma_i$s, $1 \leq i \leq p-1$ are smaller than or equal to $1/2$  there will be additional power contributions in or in front of the exponential (in the exponential the powers of the additional  terms will be smaller  than $1$). In these cases it seems difficult to get a general user-friendly formula and we suggest doing calculations on a case-by-case basis. 
\end{itemize}

\bigskip

\textbf{Case 2.} Assume that $\pi$ is \emph{infinite} such that $\pi(u,\infty)=\sum_{i=0}^p c_i u^{-\gamma_i}+O(u^{-\gamma_0+1+\varepsilon})$ as $u\rightarrow 0$, where $-1<\gamma_0-1 =\gamma_p<\gamma_{p-1}<\ldots<\gamma_1<\gamma_0<1$, $\varepsilon>0$, $c_0=1/\Gamma(1-\gamma_0)$ (this choice of $c_0$ is to simplify the expressions of the tail and density) and the others $c_i$s may be any real numbers. Then, if moreover $\gamma_1<\gamma_0-1/2$,
\begin{eqnarray*}
\mathbb P(I>t) & \underset{t \rightarrow \infty }{\sim} & k(t)t^{-\frac{\gamma_0}{1-\gamma_0}}\\
& \underset{t \rightarrow \infty }{\propto}  & t^{\frac{-\gamma_0}{2(1-\gamma_0)}-\frac{c_p\Gamma(1-\gamma_p)}{1-\gamma_0}} \exp\left(-(1-\gamma_0)t^{\frac{1}{1-\gamma_0}}-\sum_{i=1}^{p-1}\frac{c_i\Gamma(1-\gamma_i)}{1-\gamma_0+\gamma_i} t^{\frac{1-\gamma_0+\gamma_i}{1-\gamma_0}} \right).
\end{eqnarray*}
If some of the $\gamma_i$s, $1 \leq i \leq p-1$ are larger than or equal to $\gamma_0-1/2$, then there will be additional power contributions in or in front of the exponential.

\bigskip

\textbf{(\textit{a,b,c})-L\'evy measures.} A special case often encountered in practice are L\'evy measures with which are proportional to
$$e^{-ax}(1-e^{-x/c})^{b-1}\mathrm dx  \quad \text{ for some } a,c>0, b>-1$$ 
(see the next section for some examples related to random walks with a barrier, coalescence and fragmentation processes). We will call them the $(a,b,c)$-L\'evy measures. When $b \neq 0$ these are special cases of the above Case 1 and Case 2. We refer to Section \ref{sec:special}, Lemma \ref{lem:detailabc}  for details which lead to explicit expressions in terms of $a,b,c$ of the equivalents of the tail and density of $I$ in such cases. We will also see there that when $a=1$ and $b=-c$, the random variable $I$ is in fact a multiple of a generalized Mittag-Leffler random variable with parameters $(c,c)$, a random variable for which a series expansion of the density is known. Another case where the distribution of such $I$ is explicit is when $b=-1+(a+1)^{-1}$ and $c=(a+1)^{-1}$, for some $a>0$. Then, $I$ is proportional to $\mathbf e(1)^{1/(a+1)}$, where $\mathbf e(1)$ is an exponential random variable with parameter $1$.

\subsection{Some examples of applications}

We start with two classical examples of subordinators, the Gamma and stable subordinators, and then turn to more general applications related to the connections between exponential functional of subordinators and non-negative self-similar Markov processes. To simplify the exposition, we will generally only give the asymptotic of the tail $\mathbb P(I>t)$, that of the density following directly from the above results. The multiplicative constant is not in general the $\mathrm c_I$ of the theorem, but a multiple of it.


\bigskip

\textbf{Gamma subordinators.} Here $\pi(\mathrm du)=u^{-1}e^{-u}\mathrm du$ and $\phi(x)=\ln(1+x)$, so we are not in the asymptotic power expansions discussed above. From Theorem \ref{thm:main} (or Corollary \ref{cor:RV}) we obtain 
\begin{eqnarray}
\label{ex:Gammasub}
\mathbb P(I>t) &~&  \underset{t \rightarrow \infty}\propto  ~ (\ln(t))^{-1/2} \exp\left(\int_2^{t} W_{-1}\left(-\frac{e^{-\frac{1}{r}}}{r} \right) \mathrm dr + \ln(t)\right) \\
\nonumber 
&~&  \underset{t \rightarrow \infty}\propto  ~ \exp\left( -t \ln(t) - t \ln (\ln (t)) +t +O\left( \frac{t \ln (\ln (t))}{\ln(t)} \right)\right)
\end{eqnarray}
where in the first line $W_{-1}:[e^{-1},0[\rightarrow ]-\infty,-1]$ is the branch with values in $]-\infty,-1]$ of the Lambert $W$ function which is the inverse of the function $xe^x$. Series expansions of this function are obtained in \cite{Wfunc}, so one can obtain an expansion at any order of the term in the exponential of (\ref{ex:Gammasub}). To illustrate, we retain here that as $x \rightarrow 0-$, $W_{-1}(x)=\ln(-x)-\ln(-\ln(-x))+O(\ln(-\ln(-x))/\ln(-x))$. To get (\ref{ex:Gammasub}) from Theorem \ref{thm:main} just note that here $\psi$ is the inverse of $x \mapsto x/\ln(1+x)$, with $x_{\psi}=1$. Starting from $\psi(t)=t \ln (1+\psi(t))$ and setting $h(t):=-(1+\psi(t))/t$, we see that the function $h$ satisfies $h(t)e^{h(t)}=-e^{-1/t}/t$ for $t$ large enough, and that $h(t) \rightarrow -\infty$ as $t \rightarrow \infty$. Hence $h(t)=W_{-1}(-e^{-1/t}/t)$ and $$\psi(t)=-t W_{-1}(-e^{-1/t}/t) -1=t \ln (t)+ t \ln (\ln (t)) + O \left(\frac{t \ln (\ln (t))}{\ln(t)} \right) \quad \text{ as } t \rightarrow \infty,$$
which gives (\ref{ex:Gammasub}) and the line below.

\bigskip

\textbf{Stable subordinators.} When $\xi$ is a stable subordinator with Laplace exponent $\phi(x)=x^{\alpha}$, $\alpha \in (0,1)$, we immediately get from Corollary \ref{cor:RV} that as  $t \rightarrow \infty$
$$
\mathbb P(I>t) ~ =  ~ \mathrm c ~ t^{-\frac{\alpha}{2(1-\alpha)}} \exp\Big(-(1-\alpha)t^{\frac{1}{1-\alpha}} \Big) \Big(1+O\left(t^{-\frac{1}{1-\alpha}}\right) \Big)$$ 
for some $\mathrm c \in (0,\infty)$.

\bigskip

\textbf{Extinction time and moments of a non-increasing, non-negative self-similar Markov process.} Fix $\alpha>0$ and let $X:[0,\infty)\rightarrow [0,\infty)$ be a c\`adl\`ag strong Markov process which is  $1/\alpha-$self-similar, i.e. such that for all $c>0$
$$
\big\{X^{(c)}(t),t \geq 0\big\} ~\text{ is distributed as } ~\big\{cX^{(1)}(c^{-\alpha}t),t \geq 0 \big\} 
$$
where $X^{(c)}$ denotes a version of $X$ started at $c$.
Assume moreover that $X$ is non-increasing and absorbed at 0. Then it is well-known, by Lamperti transformation \cite{Lamp72}, that $X^{(1)}=\exp(-\xi_{\rho(\cdot)})$ for some subordinator $\xi$ where $\rho$ is the right-inverse of $t\mapsto \int_0^t \exp(-\alpha\xi_r) \mathrm d r$. Consequently, $X^{(1)}$ reaches 0 at time $$I_X:=\int_0^{\infty}\exp(-\alpha \xi_r) \mathrm dr,$$ of which we have an equivalent of the tail and density when $\xi$ is drift-free and unkilled, with a Laplace exponent $\phi$ satisfying (\ref{hyp:main}) (this is Theorem \ref{thm:main}). Together with Theorem 1.2 of \cite{HR12} which describes the behavior as $t \rightarrow \infty$ of $X(t)$ conditioned on being positive, this gives (see Section \ref{sec:ssMp} for details):
\begin{cor}
\label{cor:ssMp}
Assume \emph{(\ref{hyp:main})}. Then for all $a>0$,
$$
\mathbb E\big[\big(X^{(1)}(t)\big)^a\big] ~ \underset{t \rightarrow \infty}\propto ~\left(\frac{t}{\psi(\alpha t)}\right)^{1+\frac{a}{\alpha}}(\psi'(\alpha t))^{1/2}~ \exp\left(- \frac{1}{\alpha}\int_{x_{\psi}+1}^{\alpha t} \frac{\psi(r)}{r} \mathrm dr\right),
$$
where $\psi$ is the inverse of $x \mapsto x/\phi(x)$.
\end{cor}

\bigskip

\textbf{Extinction times of typical fragments in self-similar fragmentation processes.}  Some exponential functionals of subordinators appear as extinction times of typical fragments in self-similar fragmentations with a negative index, whose temporal evolution is described by a non-increasing, non-negative self-similar Markov process. See Bertoin \cite{BertoinSSF}, and also \cite{HM04} since in a related way these random variables are also the heights of uniform leaves in self-similar fragmentation trees, the first examples of which are the Brownian tree of Aldous \cite{Ald91a,Ald93} and the stable L\'evy trees of Duquesne and Le Gall \cite{DLG02,DLG05}. We will not go into details here and refer to the paper \cite{H20b} in preparation which studies the tail of the extinction time of such a fragmentation process and its connection with the extinction time of a typical fragment, and so relies partly on the present Theorem \ref{thm:main}. Specific examples will be treated there. Let us however highlight here that in several examples of fragmentation processes, the L\'evy measure associated to the extinction time of a typical fragment / height of a typical leaf belongs to the $(a,b,c)$-L\'evy measures family for which we are able to get explicit asymptotics for the tail and density: this is for example the case for Ford's Alpha-model or Aldous's Beta model (see e.g. \cite{Ford05,Ald96, HMPW08} for definitions) or for the 
Brownian and stable L\'evy trees (where, in fact, we are in the more specific case of a generalized Mittag-Leffler distribution). Again, we refer to \cite{H20b} for details.

\bigskip

\medskip

Non-negative self-similar Markov process arise naturally as scaling limits of rescaled discrete Markov chains, and relatively their extinction times describe the scaling limits of absorption times of some discrete Markov chains. We discuss here two examples related to the general case of non-increasing, $\mathbb Z_+$-valued Markov chains studied in \cite{HM11}.

\bigskip

\textbf{Number of collisions in a $\Lambda-$coalescent.} Exponential functionals of subordinators appear as scaling limits of the number of collisions in $\Lambda-$coalescent.
We briefly recall that these are Markovian models for the coalescence of particules where at each step, starting from say $n$ particles, $n-k+1$ merge into one (the other remain constant) to give $k$ distinct particules at rate
$$
\binom{n}{k-1}\int_{[0,1]} x^{n-k-1}(1-x)^{k-1} \Lambda(\mathrm dx)
$$
where $\Lambda$ is the parameter of the model, a finite measure on $[0,1]$. These models where introduced by \cite{PitCoag99,Sag99} and deeply studied since then. Let $I_n$ denote the number of steps to go from $n$ to 1 particules and set $s(u)=\int_u^1 x^{-2} \Lambda(\mathrm dx)$, $u \in (0,1)$. From \cite{IM08,HM11} we know that when $s$ is regularly varying at 0 with index $-\gamma \in (-1,0)$ and $\int_0^1 x^{-1} \Lambda(\mathrm dx)<\infty$, 
$$
\frac{I_n}{s(1/n)} \underset{n \rightarrow \infty}{\overset{\mathrm{(d)}}\longrightarrow} I
$$
where $I$ is the exponential functional of a subordinator with Laplace exponent $$\phi(x)=\frac{1}{\Gamma(2-\gamma)} \int_0^1\big(1-(1-x)^{\gamma \lambda}\big)x^{-2} \Lambda(\mathrm dx).$$ We can therefore obtain from Corollary \ref{cor:RV} the asymptotic behavior of the tail and density of this limit according to the behavior of $s$ at 0. More specifically, for the  Beta$(\alpha,\beta)-$coalescent (i.e. when the measure $\Lambda$ is the Beta$(\alpha,\beta)$ distribution, with $1<\alpha<2, \beta>0$, the number of collisions $I_n$ divided by $n^{2-\alpha}$ converge in distribution towards an  exponential functional of a subordinator with L\'evy measure
$$
\pi(\mathrm dx)=\frac{\mathrm dx}{\Gamma(a)} e^{-\beta x/(2-a)}\big(1- e^{-x/(2-a)}\big)^{3-a},
$$
so we are in the $(a,b,c)$-L\'evy measure case, with $a=\beta/(2-a)$ and $b=a-2=-c$. Lemma \ref{lem:detailabc} (in Section \ref{sec:special}) gives the behavior of the tail and density of the limiting exponential functional when $b<-1/2$ (equivalently $\alpha<3/2$) and we can obtain further results for larger $\alpha$ by using Lemma \ref{lem:abBeta} and adapting the proof of Lemma \ref{prop:aenuinfinite}.

\bigskip

\textbf{Random walks with a barrier.} To illustrate the importance of the $(a,b,c)$-L\'evy family, we finish this series of examples by recalling briefly that exponential functionals with such a L\'evy measure also arise in the context of random walks with a barrier, which are random walks conditioned on the event that they will not exceed a given threshold $n$. Several conditionings are possible, we refer to \cite[Section 2.1]{HM11} for details. Fix such a walk and assume that it is non-decreasing, with a step distribution that varies regularly with index $-c \in (-1,0)$. Then letting $n \rightarrow \infty$ (\cite{IM08,HM11})  the absorbing time of such the walk appropriately rescaled converges in distribution to the exponential functional of a subordinator with a $(c^{-1},-c,c)$-L\'evy measure with some $c \in (0,1)$, for which, again we can compute explicitly the asymptotic equivalents of the tail and density. 

\section{Preliminaries: some basic properties of $\phi$ and $\psi$}
\label{sec:phi,psi}

We gather here some properties of $\phi$ and $\psi$ that we will need in the following. We recall that as a L\'evy measure of a subordinator, $\pi$ is defined and sigma-finite on $(0,\infty)$, and such that $\int_0^{\infty}(1 \wedge v)\pi(\mathrm d v)<\infty$. Throughout we assume that $\pi(0,\infty)>0$. From its definition (\ref{def:phi}), $\phi$ is continuous increasing on $[0,\infty)$ and infinitely differentiable on $(0,\infty)$, with its $i-$th derivative given by
$$
\phi^{(i)}(x)=(-1)^{i-1}\int_0^{\infty} \exp(-xv) v^{i}  \pi(\mathrm dv), \quad i \geq 1.
$$
As immediate consequences, we have that:
\begin{enumerate}
\item[$\bullet$] $\phi$ is increasing and concave on $[0,\infty)$, the function $x \mapsto \phi(x)/x$ is decreasing on $(0,\infty)$ and converges to 0 as $x \rightarrow \infty$, and $0 \leq \phi'(x) \leq \phi(x)/x$, $\forall x>0$  
\item[$\bullet$] $\psi$, the inverse of the increasing function $x \mapsto x/\phi(x)$, is defined and increasing on $(x_{\psi},\infty)$, where $x_{\psi}:=\lim_{x \downarrow 0} x/\phi(x)=\left(\int_0^{\infty} v \pi(\mathrm dv)\right)^{-1} \in [0,\infty)$; moreover $\psi(x) \rightarrow \infty$ as $x \rightarrow \infty$
\item[$\bullet$] $\psi$ is infinitely differentiable on $(x_{\psi},\infty)$
\item[$\bullet$] $x \mapsto \psi(x)/x=\phi(\psi(x))$ is increasing on $(x_{\psi},\infty)$; in particular $x=O(\psi(x))$ as $x \rightarrow \infty$.
\end{enumerate}

\bigskip

We will need the following identities and inequalities.

\bigskip

\begin{lemma} 
\label{prop:derivative}
\begin{enumerate}
\item[\emph{(i)}] For $x>x_{\psi}$,
\begin{equation*}
\int_0^{\infty} \exp(-\psi(x)v)v \pi(\mathrm dv)=\frac{1}{x} \left(1- \frac{\psi(x)}{x\psi'(x)}\right).
\end{equation*}
\item[\emph{(ii)}] For $x>x_{\psi}$,
\begin{equation*}
\int_0^{\infty} \exp(-\psi(x)v)v^2 \pi(\mathrm dv)= \frac{2}{x^2 \psi'(x)} \left(1- \frac{\psi(x)}{x\psi'(x)}\right)-\frac{\psi(x)\psi''(x)}{x^2(\psi'(x))^3}.
\end{equation*}
\item[\emph{(iii)}] For all $a \in [1,\infty)$, 
$$
\int_0^{\infty} \exp(-\psi(x)v) v^a \pi(\mathrm dv) \lesssim \frac{1}{ x \psi(x)^{a-1}} \quad \text{on } (x_{\psi}, \infty)
$$
\item[\emph{(iv)}] For all $a \in [0,\infty)$ and all $\varepsilon \in (0,1)$
$$
\int_1^{\infty} \exp(-\psi(x)v) \exp(av) \pi(\mathrm dv) =O \left(\exp(-(1-\varepsilon) \psi(x)) \right)  \quad \text{as } x \rightarrow \infty.
$$
\end{enumerate}
\end{lemma}

\bigskip

\textbf{Proof.} We work on $(x_{\psi}, \infty)$.

\smallskip

(i) It suffices to differentiate the identity $\psi(x)=x \phi(\psi(x))$ to get:
$$
\psi'(x)=\frac{\psi(x)}{x}+x\psi'(x)\int_0^{\infty} \exp(-\psi(x)v)v \pi(\mathrm dv).
$$

(ii) Differentiating a second time,
\begin{eqnarray*}
\psi''(x)&=&\frac{\psi'(x)}{x}-\frac{\psi(x)}{x^2}+\left(\psi'(x)+x\psi''(x)\right)\int_0^{\infty} \exp(-\psi(x)v)v \pi(\mathrm dv)\\ &&-~x (\psi'(x))^2 \int_0^{\infty} \exp(-\psi(x)v)v^2 \pi(\mathrm dv)\\
&\underset{\text{by (i)}}=&  \frac{2 \psi'(x)}{x} \left(1- \frac{\psi(x)}{x\psi'(x)}\right)+\psi''(x)\left(1-\frac{\psi(x)}{x\psi'(x)}\right) -x (\psi'(x))^2 \int_0^{\infty} \exp(-\psi(x)v)v^2 \pi(\mathrm dv),
\end{eqnarray*}
which leads to the expected equality.

\smallskip

(iii) Fix $a \in [1,\infty)$ and use that $\exp(-y) \lesssim y^{-(a-1)}$ on $(0,\infty)$, to get that for $x \in (x_{\psi}, \infty)$
\begin{equation*}
\int_0^{\infty} \exp(-\psi(x)v)v^a \pi(\mathrm dv) \lesssim \frac{1}{\psi(x)^{a-1}} \int_0^{\infty} \exp\left(-\frac{\psi(x)}{2}v\right)v \pi(\mathrm dv) =
\frac{1}{\psi(x)^{a-1}}  \phi'\left(\frac{\psi(x)}{2}\right).
 \end{equation*}
Then use that $\phi'(y) \leq \phi(y)/y$ for $y>0$, that $\phi$ is increasing and the definition of $\psi$ to get
$$
 \phi'\left(\frac{\psi(x)}{2}\right) \leq  \frac{2}{\psi(x)}\phi\left(\frac{\psi(x)}{2}\right) \leq  \frac{2\phi(\psi(x))}{\psi(x)}=\frac{2}{x}.
$$ 

(iv) Fix $a \in (0,\infty), \varepsilon \in (0,1)$. Since $\psi(x) \rightarrow \infty$ as $x \rightarrow \infty$, ~$\psi(x)\geq a/\varepsilon$~ for all $x$ large enough. Thus, for such $x$s and all $v\geq 0$, $\exp(-\psi(x)v) \exp(av) \leq \exp(-(1-\varepsilon)\psi(x)v)$, and then
$$
\int_1^{\infty} \exp(-\psi(x)v) \exp(av) \pi(\mathrm dv) \leq \pi(1,\infty)\exp(-(1-\varepsilon)\psi(x)).
$$
$\hfill \square$

\bigskip

\begin{lemma}
\label{prop:contract}
Assume that $
\frac{x\phi'(x)}{\phi(x)} \leq c <1
$
for $x$ large enough. Then for $x$ large enough
$$
\int_0^{\infty} \exp(-xv)v \pi(\mathrm dv) \leq \frac{c \phi(x)}{x}
$$
or equivalently, using that $\phi(\psi(x))/\psi(x))=x^{-1}$,
$$
\int_0^{\infty} \exp(-\psi(x)v)v \pi(\mathrm dv) \leq \frac{c}{x}.
$$
\end{lemma}

\textbf{Proof.} Just use that 
$$\int_0^{\infty} \exp(-xv)v \pi(\mathrm dv)=\phi'(x) \leq c \frac{\phi(x)}{x}$$ 
for all $x$ large enough. $\hfill \square$

\bigskip

Next, assuming (\ref{hyp:main}) only, we can control the derivatives of $\psi$ as follows.

\begin{lemma}
\label{prop:boundpsi}
Assume \emph{(\ref{hyp:main})}. Then as $x \rightarrow \infty$
$$
 \psi'(x) =O\left(\frac{\psi(x)}{x}\right) \quad \text{and} \quad \left|\psi''(x) \right| =O\left( \frac{\psi(x)}{x^2}\right).
$$
Consequently, there exists $\mathrm c \in (0,\infty)$ such that for all $x$ large enough and all $v \geq 0$,
$\psi(xe^v) \leq e^{\mathrm cv} \psi(x).$
\end{lemma}

(We recall that $\psi'$ is positive since $\psi$ is increasing.)

\bigskip

\textbf{Proof.} From the equality $\psi(x)=x\phi(\psi(x))$, $x>x_{\psi}$, we get by differentiating twice 
\begin{equation}
\label{eq:hyp}
\psi'(x)=\phi(\psi(x))+x \phi'(\psi(x))\psi'(x) 
\end{equation}
\begin{equation*}
\psi''(x)= 2\phi'(\psi(x)) \psi'(x) +  x \phi''(\psi(x)) (\psi'(x))^2+x\phi'(\psi(x)) \psi''(x).
\end{equation*}
Besides, from (\ref{hyp:main}), there exists $c_{\phi} \in (0,1)$ and $x_{\phi} \in (0,\infty)$ such that $\phi'(x) \leq c_{\phi} \phi(x)/x$ for all $x >x_{\phi}$, hence $x\phi'(\psi(x)) \leq  c_{\phi} x \phi(\psi(x))/\psi(x)=  c_{\phi}$ for all $x>x_{\phi}/\phi(x_{\phi}) \geq x_{\psi}$. 

\smallskip

(i) We therefore get that for  $x>x_{\phi}/\phi(x_{\phi})$
\begin{eqnarray*}
\psi'(x) \leq  \phi(\psi(x))+ \mathrm{c_{\phi}} \psi'(x)
\end{eqnarray*}
and then
\begin{equation}
\label{eq:psiprime}
\psi'(x) \leq \frac{1}{1-\mathrm c_{\phi}} \cdot \frac{\psi(x)}{x} .
\end{equation}
This implies that for $v \geq 0$ and $x>x_{\phi}/\phi(x_{\phi})$,
$$
\ln \left(\frac{\psi(xe^v)}{\psi(x)}\right) = \int_x^{xe^v} \frac{\psi'(u)}{\psi(u)} \mathrm du  \leq \frac{1}{1-\mathrm c_{\phi}} \int_x^{xe^v} \frac{1}{u} \mathrm du = \frac{v}{1-\mathrm c_{\phi}},
$$
hence $\psi(xe^v) \leq e^{ \frac{v}{1-\mathrm c_{\phi}}} \psi(x)$.

\smallskip

(ii) Next to get the bound for $|\psi''|$, we use the above equation involving $\psi''$ to get for $x>x_{\psi}$
$$
\left| \psi''(x)\right| \cdot \left| 1-x\phi'(\psi(x))\right| \leq 2\phi'(\psi(x)) \psi'(x) + x \left|\phi''(\psi(x))\right| (\psi'(x))^2.
$$
Since $ x\phi'(\psi(x)) \leq \mathrm c_{\phi}<1$ for $x>x_{\phi}/\phi(x_{\phi})$ and  $|\phi''(\psi(x))\big| \lesssim 1/x\psi(x)$ for $x>x_{\psi}$ by Lemma \ref{prop:derivative} (iii) (applied to $a=2$), this implies that for $x>x_{\phi}/\phi(x_{\phi})$
\begin{eqnarray*}
\left|\psi''(x) \right| \lesssim  \frac{\psi'(x)}{x}+ \frac{(\psi'(x))^2}{\psi(x)} \lesssim \frac{\psi(x)}{x^2},
\end{eqnarray*}
where we used (\ref{eq:psiprime}) for the second inequality.
$\hfill \square$

\bigskip

Last, we note that in the regular variation case we have:

\begin{lemma}
\label{cororegvar}
If $\phi$ is regularly varying at $\infty$ with index $\gamma \in [0,1)$, then $\psi$ is regularly varying at $\infty$ with index $(1-\gamma)^{-1}$ and 
$$
\frac{x\psi'(x)}{\psi(x)} \underset{x \rightarrow \infty}\longrightarrow \frac{1}{1-\gamma}.
$$
\end{lemma}

\textbf{Proof.} Since  $\phi$  is regularly varying at $\infty$ with index $\gamma \in [0,1)$ and $\phi'$ is decreasing, we know by the Monotone Density Theorem (\cite[Theorem 1.7.2]{BGT}) that 
$$
\frac{x\phi'(x)}{\phi(x)}  \underset{x \rightarrow \infty}\longrightarrow \gamma.
$$
Then we conclude with (\ref{eq:hyp}),
$$
 \left(1- \frac{\psi(x)}{x\psi'(x)}\right)= \frac{\psi(x)\phi'(\psi(x))}{\phi(\psi(x))}  \underset{x \rightarrow \infty}\longrightarrow \gamma.
$$
$\hfill \square$

\section{Proof of Theorem \ref{thm:main}}
\label{sec:thmmain}

Let $$f(x):=-\ln \mathbb P(I>x), \quad x\geq 0.$$ 
By Carmona, Petit and Yor \cite[Proposition 2.1.]{CPY97} and Pardo, Rivero and Schaik \cite[Theorem 2.3]{PRvS13}, we know that the law of $I$ has a density, denoted by $k$, with respect to the Lebesgue measure on $[0,\infty)$ which satisfies the integral equation
\begin{equation}
\label{eq:k}
k(x)=\int_x^{\infty} \bar \pi\left(\ln(y/x)\right) k(y) \mathrm dy = \int_0^{\infty} \left(\int_x^{xe^v} k(y) \mathrm dy \right) \pi(\mathrm dv)
\end{equation}
where  $\bar \pi(z)=\int_z^{\infty} \pi(\mathrm dv), z>0$. Morever $k$ is the unique non-negative solution to its equation such that $\int_0^{\infty}k(x)\mathrm dx=1$. Patie and Savov \cite[Theorem 2.4]{PS18} proved that $k$ is infinitely differentiable on $(0,\infty)$, at least under our zero-drift assumption. Consequently, $f$ is infinitely differentiable on $(0,\infty)$ and its derivative $f'(x)=k(x)/\mathbb P(I>x)$, $x>0$, satisfies the integral equation
\begin{equation}
\label{eqfprime}
f'(x)=\int_0^{\infty} \left(1-\exp\left(-\int_x^{xe^v}f'(u)\mathrm du \right) \right) \pi(\mathrm dv), \quad \forall x >0.
\end{equation}
Reasoning informally, one sees that for small $v$s the integral $\int_x^{xe^v}f'(u)\mathrm du$ is close to $f'(x) x v $, so in some rough sense the equation (\ref{eqfprime}) satisfied by $f'$ is close  to the one satisfied by the function $x \mapsto \psi(x)/x$ on its domain of definition:
$$
\frac{\psi(x)}{x}=\phi(\psi(x))=\int_0^{\infty} \left(1-\exp\left(-\frac{\psi(x)}{x}xv\right)\right) \pi(\mathrm dv), \quad \forall x >x_{\psi}.
$$
This ``proximity" of $f'(x)$ and $\psi(x)/x$ for large $x$ was used in \cite{HR12} to get the extension (\ref{extR}) of Rivero's result. Our key result in the present paper is to get the following precise asymptotic expansion of $f'$ in terms of $\psi$:

\begin{prop}
\label{prop:equiv}
Assume \emph{(\ref{hyp:main})}.
Then as $x \rightarrow \infty$
$$f'(x)= \frac{\psi(x)}{x}+\frac{\psi'(x)}{\psi(x)} -\frac{1}{x} - \frac{\psi''(x)}{2 \psi'(x)} +O\left(\frac{\psi'(x)}{ (\psi( x))^2}\right).$$
\end{prop}

\bigskip

Note that this proposition immediately implies Theorem \ref{thm:main}. Indeed, 
$$
\int_{x_{\psi+1}}^{\infty} \left | \frac{\psi'(x)}{ (\psi( x))^2}\right | \mathrm dx= \int_1^{\infty}  \frac{\psi'(x)}{ (\psi( x))^2}  \mathrm dx= \frac{1}{\psi(1)}
$$
since $\psi'$ is positive and $\psi(x) \rightarrow \infty$ as $x \rightarrow \infty$. By Proposition \ref{prop:equiv} we therefore have 
\begin{eqnarray*}
 -\ln \mathbb P(I>t)=f(t)=\int_{x_{\psi}+1}^t \frac{\psi(x)}{x} \mathrm dx+\ln({\psi(t)})-\ln(t)-\frac{1}{2}\ln (\psi'(t))+O(1) \quad \text{as } t \rightarrow \infty,
\end{eqnarray*}
which gives Theorem \ref{thm:main}. 

\bigskip

The rest of this section is devoted to the proof of Proposition \ref{prop:equiv}. In that aim we introduce the function 
$$
h(x):=f'(x)-\frac{\psi(x)}{x}, 
$$ 
which is defined and infinitely differentiable on $(x_{\psi},\infty)$.  Our proof consists first of setting up some rough power bounds for $h,h'$ as $x \rightarrow \infty$ (Section \ref{sec:hbounds}) and useful estimates (Section \ref{sec:estimates}). Then use both of them to prove that $h(x)=O(x^{-1})$ and $h'(x)=O(x^{-2})$ as $x \rightarrow \infty$ (Section \ref{sec:refin}), which are the right-order bounds for $h$ and $h'$. And finally use these new bounds to refine the estimates and  improve the asymptotic expansion of $h$ to get Proposition \ref{prop:equiv} (Section \ref{sec:proofProp}).

\bigskip

Throughout the rest of this section we assume that Hypothesis (\ref{hyp:main}) holds and fix once and for all $\kappa \in (0,1)$ such that 
\begin{equation}
\label{def:kappa}
\frac{x\phi'(x)}{\phi(x)} < \kappa
\end{equation}
for all $x$ large enough. 

\subsection{First bounds for $h, h'$}
\label{sec:hbounds}

The goal of this section is to establish power uper bounds for $|h|,|h'|$. These first bounds are rough but needed to establish the finer bounds we will need for the proof of Proposition \ref{prop:equiv}. So, recalling the above definition of $\kappa$, we will show here that

\begin{prop}
\label{prop_majo_hprime}
$h(x) =O\big( x^{\frac{\kappa}{1-\kappa}}\big)$ and $h'(x) =O\big( x^{\frac{2\kappa}{1-\kappa}}\big)$ as $x \rightarrow \infty$.
\end{prop}

We will also need the positivity of $h$ in a neighborhood of $+\infty$:

\begin{prop}
\label{lem:hpositive}
There exists $x_{+} \in (x_{\psi},\infty)$ such that $h$ is positive on $[x_{+},\infty)$. 
\end{prop}

To set up these results, we will use some sequences of functions built iteratively that converge respectively to $f'$ and $x \mapsto \psi(x)/x.$ This approach was already used in \cite{HR12} and in fact the positivity of $h$ and the fact that $h(x) =O\big( x^{\frac{\kappa}{1-\kappa}}\big)$  as $x \rightarrow \infty$ can be read in the proofs there.  So the only real new result here is that $h'(x) =O\big( x^{\frac{2\kappa}{1-\kappa}}\big)$ as $x \rightarrow \infty$. We give however complete proofs of Proposition \ref{prop_majo_hprime} and Proposition \ref{lem:hpositive}, to ease the reading and since we use a slightly different and simpler approach.

\subsubsection{Approximation schemes for $f'$ and $\psi(x)/x$} 

\textbf{Two operators.} On the set of non-negative continuous functions  defined on $[0,\infty)$ we consider two operators $\Theta$ and $\Theta_{\phi}$ defined for any non-negative continuous function $g$ by
$$
\Theta(g)(x):=\int_0^{\infty} \left(1-\exp\left(-\int_x^{xe^v} g(u) \mathrm du\right) \right) \pi(\mathrm dv), \quad \forall x \geq 0
$$
(note that $\Theta(g)(x) \in [0,\infty)$ for all $x\geq 0$)
and
$$
\Theta_{\phi}(g)(x):=\phi\left(xg(x)\right), \quad \forall x \geq 0.
$$
We note that $\Theta(g)$ and $\Theta_{\phi}(g)$ are also non-negative continuous functions on $[0,\infty)$. Moreover these operators are non-decreasing: if $g_1(x) \leq g_2(x)$ for all $x\in [0,\infty)$ then $\Theta(g_1)(x)\leq \Theta(g_2)(x)$ and $\Theta_{\phi}(g_1)(x)\leq \Theta_{\phi}(g_2)(x)$  for all $x\in [0,\infty)$.

\medskip

The function $f'$ is a fixed point of $\Theta$, whereas $x \mapsto \psi(x)/x$ is a fixed point of $\Theta_{\phi}$ restricted to functions defined on $(x_{\psi},\infty)$. Recalling the definition of $\psi$ as the inverse of $x \mapsto x/\phi(x)$, we see that $x \mapsto \psi(x)/x$ is the unique fixed point of  $\Theta_{\phi}$. The function $f'$ is also the unique fixed point of $\Theta$: let $p:[0,\infty) \rightarrow [0,\infty)$ be some continuous function such that $\Theta(p)=p$ and let $P$ denotes its primitive null at 0; then the function $r$ defined by $r(x):=p(x)\exp(-P(x)), x \geq 0$ is continuous non-negative and such that $r(x)=\int_0^{\infty} \left( \int_x^{xe^v} r(y) \mathrm dy \right) \pi(\mathrm dv) $ and $\int_0^{\infty} r(x) \mathrm dx=1$. Hence $r=k$, by uniqueness of the non-negative solutions to (\ref{eq:k}).

\bigskip

\textbf{Approximation schemes.} Now consider $a \in (0,\pi[0,\infty))$ and let $g_0(x)=g_0^{\phi}(x):=a$, $\forall x \geq 0$. Then define recursively for $n\geq 1$
$$
g_n(x):=\Theta(g_{n-1})(x), \quad x \geq 0
$$
and
$$
g^{\phi}_n(x):=\Theta_{\phi}(g^{\phi}_{n-1})(x), \quad x \geq 0.
$$
Clearly, reasoning recursively, the functions $g_n,g^{\phi}_n$ are all non-decreasing on $[0,\infty)$.
Note also that $g_1(x) \geq g^{\phi}_1(x)$ for all $x \geq 0$ (since $e^v-1\geq v$ for $v \geq 0$) and that for $x$ large enough, say $x\geq x_a$, $g^{\phi}_1(x) \geq g^{\phi}_0(x)$ since $\phi(xa) \rightarrow  \pi[0,\infty) > a$ as $x \rightarrow \infty$. So for all $x \in [x_a,\infty)$
$$
g_1(x) \geq g^{\phi}_1(x), \quad g_1(x) \geq g_0(x), \quad  g^{\phi}_1(x)\geq g^{\phi}_0(x).
$$
Iterating, using that both operators are non-decreasing and that the functions $g_n,g^{\phi}_n$ are all non-decreasing, we see that
$$
g_n(x) \geq g^{\phi}_n(x), \quad g_n(x) \geq g_{n-1}(x), \quad  g^{\phi}_n(x)\geq g^{\phi}_{n-1}(x), \quad \forall x \geq x_a, \forall n \geq 1.
$$

Besides, there exists $x_{\kappa}$ such that 
\begin{equation}
\label{majog_n}
g_n(x) \leq x^{\frac{\kappa}{1-\kappa}} \quad \text{ for all } x \geq x_{\kappa} \text{ and all } n\geq 1.
\end{equation}
This can be proved by induction on $n$. We start by choosing $x_{\kappa} \geq a^{1-\kappa}$ and such that 
$$
\pi(1,\infty)+\phi(x^\frac{1}{1-\kappa}) + x^\frac{1}{1-\kappa}\int_0^1 \exp(-x^\frac{1}{1-\kappa}v)\left((1-\kappa) (e^{\frac{v}{1-\kappa}}-1)-v\right) \pi(\mathrm dv) \leq x^{\frac{\kappa}{1-\kappa}} \quad \forall x \geq x_{\kappa}.
$$
Such $x_{\kappa}$ exists because: 1)  by (\ref{def:kappa}) there exists $\varepsilon \in (0,\kappa)$ such that $\phi(x) =O\left(x^{\kappa-\varepsilon}\right)$ as $x \rightarrow \infty$ and 2) by Lemma \ref{prop:derivative} (iii) $\int_0^1 \exp(-xv)v^2 \pi(\mathrm dv)=O(\phi(x)/x^2)$ as $x \rightarrow \infty$ (recalling that $\psi^{-1}(x)=x/\phi(x)$) and therefore
\begin{eqnarray*}
&&\pi(1,\infty)+\phi(x^\frac{1}{1-\kappa}) +  x^{\frac{1}{1-\kappa}} \int_0^1 \exp(-x^\frac{1}{1-\kappa}v)\left((1-\kappa) (e^{\frac{v}{1-\kappa}}-1)-v\right) \pi(\mathrm dv) \\
&=&O\left(x^{\kappa-\varepsilon}\right)+x^\frac{1}{1-\kappa}O(\phi(x^\frac{1}{1-\kappa})/x^\frac{2}{1-\kappa})\\
&=&O(x^{\kappa-\varepsilon}).
\end{eqnarray*}
The proof by induction is then easy. First, $g_0(x)=a \leq x^{\frac{1}{1-\kappa}}$ for $x \geq x_{\kappa}$. Then assume that $g_n(x) \leq x^{\frac{\kappa}{1-\kappa}}$ for all $x \geq x_{\kappa}$ and note that since $\int_x^{xe^{v}} u^{\frac{\kappa}{1-\kappa}} \mathrm du = (1-\kappa)x^{\frac{1}{1-\kappa}} (e^{\frac{v}{1-\kappa}}-1)$ then for all $x \geq 0$
$$
g_{n+1}(x) \leq \pi(1,\infty)+\phi(x^\frac{1}{1-\kappa}) + x^\frac{1}{1-\kappa}\int_0^1 \exp(-x^\frac{1}{1-\kappa}v)\left((1-\kappa) (e^{\frac{v}{1-\kappa}}-1)-v\right) \pi(\mathrm dv) 
$$
which is smaller than $x^{\frac{\kappa}{1-\kappa}}$ for $x \geq x_{\kappa}$.

\bigskip

\textbf{First conclusions.} So for all $x \geq \max(x_{\kappa},x_a)$ the sequences $(g_n(x))_n$  and $(g^{\phi}_n(x))_n$ are non-decreasing and bounded. Their limits are fixed points
respectively of $\Theta$ and $\Theta_{\phi}$, hence are respectively equal to $f'$ and $x \mapsto \psi(x)/x$ (on $(x_{\psi},\infty)$).  Using the links stated above between the functions $g_n,g_n^{\phi}$ and $x \mapsto  x^{\frac{\kappa}{1-\kappa}}$, we therefore see that for all $x$ sufficiently large, 
\begin{equation}
\label{eq:boundfprime}
\frac{\psi(x)}{x}  ~ \leq ~ f'(x) ~ \leq ~ x^{\frac{\kappa}{1-\kappa}}.
\end{equation}
This immediately implies Proposition \ref{lem:hpositive} and the first assertion of Proposition \ref{prop_majo_hprime}.

\subsubsection{Upper bound for $|h'|$}

It remains to prove the second assertion of Proposition \ref{prop_majo_hprime}. The main step is to settle a similar upper bound for $f''$. In that aim we first set up two lemmas. We use the notation introduced in the approximation schemes of the previous section: $(g_n^{\phi}), (g_n), a, x_a$.

\begin{lemma}
\label{lem:bornegn}
There exists $K \in (0,\infty)$ and $x_1 \geq 0$ such that for all $x \geq x_1$ and all $n\geq 1$, 
$$
\frac{g^{\phi}_{n+1}(x)}{g^{\phi}_n(x)} \leq (Kx^{\kappa})^{\kappa^n}.
$$
\end{lemma}

\textbf{Proof.} By (\ref{def:kappa}), $\frac{\phi(y)}{\phi(x)} \leq \left(\frac{y}{x}\right)^{\kappa}$ for all $y \geq x$ large enough, say $y\geq x \geq x_0>0$. In particular $\phi(x) \leq C x^{\kappa}, \forall x \geq x_0$, with $C:=\phi(x_0)x_0^{-\kappa}$. Then we note that for $x \geq x_1:=\max(x_0/a,x_{a})$ and all $n \geq 1$
$$
\frac{g^{\phi}_{n+1}(x)}{g^{\phi}_n(x)}=\frac{\phi(xg^{\phi}_n(x))}{\phi(xg^{\phi}_{n-1}(x))}=\frac{\phi\left(x\frac{g^{\phi}_n(x)}{g^{\phi}_{n-1}(x)}g^{\phi}_{n-1}(x)\right)}{\phi(xg^{\phi}_{n-1}(x))} \leq \left( \frac{g^{\phi}_n(x)}{g^{\phi}_{n-1}(x)}\right)^{\kappa}.
$$
Iterating we get that for $x \geq x_1$
$$
\frac{g^{\phi}_{n+1}(x)}{g^{\phi}_n(x)}  \leq \left( \frac{g^{\phi}_1(x)}{g^{\phi}_{0}(x)}\right)^{\kappa^n} = \left( \frac{\phi(ax)}{a}\right)^{\kappa^n} \leq (Ca^{\kappa-1} x^{\kappa})^{\kappa^n}
$$
as expected.
$\hfill \square$

\bigskip

Now introduce for $n\geq 1$, $x \geq x_a$, $\varepsilon \geq 0$, $$F_n(x,\varepsilon):=g_n(x(1+\varepsilon))-g_n(x)$$ which is non-negative (since $g_n$ is non-decreasing on $[x_a,\infty)$).

\begin{lemma}
\label{lem1}
There exists $C \in (0,\infty)$ and $x_2 \geq 0$ such that for all $n \geq 0, x \geq x_2$ and all $\varepsilon \in [0,1]$ 
$$
0 \leq F_n(x,\varepsilon) \leq \frac{\varepsilon \left(Cx \right)^{\frac{2\kappa}{1-\kappa}}}{1-\kappa}
$$
Letting first $n \rightarrow \infty$ and then $\varepsilon \rightarrow 0$, this implies that
$$
0 ~\leq~ f''(x) ~\lesssim ~x ^{\frac{2\kappa}{1-\kappa}} \quad \text{on } [x_2,\infty).
$$
\end{lemma}

\bigskip

\textbf{Proof.} $\bullet$ We first need to fix some notation. We denote by $x'_{\kappa}$ the threshold such that 
\begin{equation}
\label{eq:vvcarre}
\int_0^{\infty} \exp(-xv)  \frac{1-\kappa}{2-\kappa}\big(e^{\frac{2-\kappa}{1-\kappa} v}-1\big) \pi(\mathrm dv) \leq \kappa \frac{\phi(x)}{x} \quad \text{ for all } x \geq x'_{\kappa}.
\end{equation}
Such $x'_{\kappa}$ exists since 
$$
\frac{1-\kappa}{2-\kappa}\big(e^{\frac{2-\kappa}{1-\kappa} v}-1\big)-v \lesssim v^2 \quad \text{on }[0,1]
$$
and by Lemma \ref{prop:contract}, Lemma \ref{prop:derivative} (iii) and (iv) (using that $\psi^{-1}(x)=x/\phi(x)$). Then, with the notation of the previous lemma, we chose $C \geq \max(K,2)$ and $x_2 \geq \max(x_{\kappa},x'_{\kappa}x_a, x_1,1)$. 

\medskip

$\bullet$ We now prove by induction on $n \geq 0$ that for all $x \geq x_2$ and all $\varepsilon \in [0,1]$ 
\begin{equation}
\label{eq:Fn}
F_n(x,\varepsilon) \leq \varepsilon\left(Cx \right)^{\frac{\kappa}{1-\kappa}}\sum_{i=1}^n \kappa^i (Cx^{\kappa})^{\sum_{k=n-i}^{n-1} \kappa^k}
\end{equation}
(by convention the sum is null when $n=0$), which clearly leads to the statement of the lemma.
The initialization is immediate since $F_0(x,\varepsilon)=0$ for all $x \geq 0$ and $\varepsilon \geq 0$. Then assume that (\ref{eq:Fn}) holds for some $n \geq 0$ and all $x \geq x_2$, $\varepsilon \geq 0$. Since $\exp(-y)-\exp(-z) \leq \exp(-y)(z-y)$ for all $z \geq y \geq 0$, we get that for $x \geq x_2$ and $\varepsilon \geq 0$,
{\setlength{\jot}{12pt}
\begin{eqnarray*}
F_{n+1}(x,\varepsilon) & = & \Theta(g_n)(x(1+\varepsilon)) - \Theta(g_n)(x) \\
&\leq& \int_0^{\infty} \exp\left(-\int_x^{xe^v} g_n(u) \mathrm du\right) \left( \int_{x(1+\varepsilon)}^{x(1+\varepsilon)e^v} g_n(u) \mathrm du-\int_x^{xe^v} g_n(u) \mathrm du\right) \pi(\mathrm dv) \\
&\leq &  \int_0^{\infty} \exp\left(- g_n(x) x v\right) \left( \int_{x(1+\varepsilon)}^{x(1+\varepsilon)e^v} g_n(u) \mathrm du-\int_x^{xe^v} g_n(u) \mathrm du\right) \pi(\mathrm dv)
\end{eqnarray*}}
where we have used that $g_n$ is non-decreasing and positive  on $[x_2,\infty)$ (which in particular implies that the difference $ \int_{x(1+\varepsilon)}^{x(1+\varepsilon)e^v} g_n(u) \mathrm du-\int_x^{xe^v} g_n(u) \mathrm du$ is positive for all $x \geq x_2,\varepsilon \geq 0$). Then by the induction hypothesis and since $g_n(x) \leq x^{\frac{\kappa}{1-\kappa}}$ ~ for $x \geq x_2$, we have for $\varepsilon \in [0,1]$
\begin{eqnarray*}
&& \int_{x(1+\varepsilon)}^{x(1+\varepsilon)e^v} g_n(u) \mathrm du -\int_x^{xe^v} g_n(u) \mathrm du =  \int_{x}^{xe^v} g_n(u(1+\varepsilon)) (1+\varepsilon)\mathrm du-\int_x^{xe^v} g_n(u) \mathrm du \\ 
& \leq& \varepsilon \left(Cx \right)^{\frac{\kappa}{1-\kappa}}\sum_{i=1}^n \kappa^i \int_{x}^{xe^v}(Cu^{\kappa})^{\sum_{k=n-i}^{n-1} \kappa^k}  \mathrm du + \varepsilon  (1+\varepsilon)^{\frac{\kappa}{1-\kappa}} \int_{x}^{xe^v} u^{\frac{\kappa}{1-\kappa}} \mathrm du \\
 & \leq & \varepsilon \left(Cx \right)^{\frac{\kappa}{1-\kappa}}\sum_{i=1}^n \kappa^i (Cx^{\kappa})^{\sum_{k=n-i}^{n-1}\kappa^k}  x \frac{e^{v\left(\sum_{k=n-i}^{n-1} \kappa^k+1\right)}-1}{\sum_{k=n-i}^{n-1} \kappa^k+1} 
 + \varepsilon   (1+\varepsilon)^{\frac{\kappa}{1-\kappa}} x^{\frac{\kappa}{1-\kappa}}  x\frac{(e^{v\left(\frac{1}{1-\kappa}\right)}-1)}{\frac{1}{1-\kappa}} \\
& \leq & \varepsilon \left(Cx \right)^{\frac{\kappa}{1-\kappa}} \left(\sum_{i=1}^n \kappa^i (Cx^{\kappa})^{\sum_{k=n-i}^{n-1}\kappa^k} +1\right)  \frac{1-\kappa}{2-\kappa}(e^{v\left(\frac{2-\kappa}{1-\kappa}\right)}-1)
 \end{eqnarray*}
 where for the last inequality we use that $C \geq 2$ and both $\frac{e^{v\left(\sum_{k=n-i}^{n-1} \kappa^k+1\right)}-1}{\sum_{k=n-i}^{n-1} \kappa^k+1} \leq \frac{1-\kappa}{2-\kappa}(e^{v\left(\frac{2-\kappa}{1-\kappa}\right)}-1)$ and $(1-\kappa)(e^{v\left(\frac{1}{1-\kappa}\right)}-1)\leq  \frac{1-\kappa}{2-\kappa}(e^{v\left(\frac{2-\kappa}{1-\kappa}\right)}-1)$ since $\sum_{k=n-i}^{n-1}\kappa^k+1 \leq \frac{2-\kappa}{1-\kappa}$ for all  $0\leq i \leq n$. 
 
 Thus we have for $x \geq x_2$,
 \begin{eqnarray*}
 && F_{n+1}(x,\varepsilon) \\
 &\leq&  \varepsilon \left(Cx \right)^{\frac{\kappa}{1-\kappa}} \left(\sum_{i=1}^n \kappa^i (Cx^{\kappa})^{\sum_{k=n-i}^{n-1}\kappa^k} +1\right) x \int_0^{\infty} \exp\left(- g_n(x) x v\right)  \frac{1-\kappa}{2-\kappa}(e^{v\left(\frac{2-\kappa}{1-\kappa}\right)}-1) \pi(\mathrm dv) \\
 &\underset{\text{by } (\ref{eq:vvcarre})}\leq&  \varepsilon \left(Cx \right)^{\frac{\kappa}{1-\kappa}} \left(\sum_{i=1}^n \kappa^i (Cx^{\kappa})^{\sum_{k=n-i}^{n-1}\kappa^k} +1\right)   \kappa  \frac{\phi(xg_n(x))}{g_n(x)}.
 \end{eqnarray*}
Since the function $x \mapsto \phi(x)/x$ is decreasing and $g^{\phi}_{n} \leq g_n$, we have for $x \geq x_2$, 
$$
\frac{\phi(xg_n(x))}{g_n(x)} \leq \frac{\phi(xg^{\phi}_n(x))}{g^{\phi}_n(x)}=\frac{g^{\phi}_{n+1}(x)}{g_n^{\phi}(x)} \leq (Cx^{\kappa})^{\kappa^n}
$$
where for the last inequality we use Lemma \ref{lem:bornegn} and that  $K \leq C$. This implies that
\begin{eqnarray*}
F_{n+1}(x,\varepsilon) &\leq& \varepsilon  \left(Cx \right)^{\frac{\kappa}{\kappa-1}}\left(\sum_{i=1}^{n} \kappa^{i+1}  (Cx^{\kappa})^{\sum_{k=n-i}^{n} \kappa^k} + \kappa(Cx^{\kappa})^{\kappa^n}\right) \\
&=&  \varepsilon  \left(Cx \right)^{\frac{\kappa}{\kappa-1}}\left(\sum_{i=1}^{n+1} \kappa^{i}  (Cx^{\kappa})^{\sum_{k=n+1-i}^{n} \kappa^k}\right)
\end{eqnarray*}
as expected. $\hfill \square$

\bigskip

\textbf{Proof of Proposition \ref{prop_majo_hprime}.}
By definition of $h$, for $x>x_{\psi}$,
$$
h'(x)=f''(x)-\frac{\psi'(x)}{x}+\frac{\psi(x)}{x^2}.
$$
Together with the previous lemma, the upper bound (\ref{eq:boundfprime}) and Lemma \ref{prop:boundpsi}, this implies that 
$$
h'(x)=O\big(x^{\frac{2\kappa}{1-\kappa}}\big) \quad \text{ as } x \rightarrow \infty.
$$
$\hfill \square$

\subsection{Intermediate estimates}
\label{sec:estimates}

We will now analyze more closely the equation satisfied by $h$ and in that aim we let for all $x>x_{\psi}$, $v \geq 0$,
$$A(x,v):=\int_x^{x e^v} \left(f'(u) -\frac{\psi(u)}{u} \right) \mathrm du=\int_x^{x e^v} h(u) \mathrm du$$
and
$$
B(x,v):=\int_x^{x e^v}\frac{\psi(u)}{u} \mathrm du-\psi(x)v
$$
and note that 
\begin{equation}
\label{hAB}
h(x) = \int_0^{\infty}\exp(-\psi(x)v) \left(1-\exp\left(-A(x,v)-B(x,v)\right) \right)\pi(\mathrm dv).
\end{equation}
The goal of this section is to set up some informations on the functions $A$ and $B$. We first note that $A(x,v) \geq 0$ for $x \geq x_+$ and all $v \geq 0$, by Proposition \ref{lem:hpositive} and that $B(x,v) \geq 0$ for all $x>x_{\psi},v \geq 0$ since the function $x \mapsto \psi(x)/x$ is increasing on $(x_{\psi},\infty)$ and $e^v-1-v \geq 0$ for all $v \geq 0$.

\bigskip

Then we introduce for $x>x_{\psi}$ the functions 
$$
h_{\frac{\kappa}{1-\kappa}}(x):= \frac{h(x)}{x^{\frac{\kappa}{1-\kappa}}} \quad \text{and} \quad h'_{\frac{2\kappa}{1-\kappa}}(x):= \frac{h'(x)}{x^{\frac{2\kappa}{1-\kappa}}} 
$$
and
$$
\bar h_{\frac{\kappa}{1-\kappa}}(x):=\sup_{u \geq x} |h_{\frac{\kappa}{1-\kappa}}(u)|  \quad \text{and} \quad
\bar h'_{\frac{2\kappa}{1-\kappa}}(x):=\sup_{u \geq x} |h'_{\frac{2\kappa}{1-\kappa}}(u)|.
$$
A key point is that these suprema are \emph{finite}, by Proposition \ref{prop_majo_hprime}. We will need the following bounds.

\bigskip

\begin{lemma}
\label{lem:propA}
\begin{enumerate}
\item[\emph{(i)}] $|A(x,v)|x^{-\frac{\kappa}{1-\kappa}} -\bar h_{\frac{\kappa}{1-\kappa}}(x) x v ~\lesssim~ \bar h_{\frac{\kappa}{1-\kappa}}(x) x v^2$ ~ for $(x,v) \in (x_{\psi}, \infty) \times [0,1]$
\item[\emph{(ii)}] $\left| \partial_x A(x,v)  \right| \leq |h(xe^v)e^v|+|h(x)|$ \hspace{0.1cm} for all $x>x_{\psi}$, $v \geq 0$.
\item[\emph{(iii)}] $\left| \partial_x A(x,v)  \right| x^{-\frac{2\kappa}{1-\kappa}}- \bar h'_{\frac{2\kappa}{1-\kappa}}(x) x v ~\lesssim ~ h(x)x^{-\frac{2\kappa}{1-\kappa}}v + \bar h'_{\frac{2\kappa}{1-\kappa}}(x) x v^2$ \hspace{0.1cm}  for $(x,v) \in (x_{\psi}, \infty) \times [0,1]$.
\end{enumerate}
\end{lemma}

(We recall that (i) means that there exists a finite positive constant $\mathrm c$ such that $|A(x,v)|x^{-\frac{\kappa}{1-\kappa}}  -\bar h_{\frac{\kappa}{1-\kappa}}(x) x v \leq \mathrm c \bar h_{\frac{\kappa}{1-\kappa}}(x) x v^2$ for all $x \in (x_{\psi},\infty)$ and all $v \in [0,1]$.)

\bigskip

\begin{lemma}
\label{lem:propB}
There exists $x_0\geq x_{\psi}$ and $\mathrm c >2$ such that 
\begin{enumerate}
\item[\emph{(i)}] 
$|B(x,v)-\psi'(x)x\frac{v^2}{2}| ~\lesssim~ \psi(x) v^3 \hspace{0.1cm} \text{ for } (x,v) \in (x_0, \infty) \times [0,1]$.

Consequently, $0 \leq B(x,v)   ~\lesssim~   \psi(x) v^2$ \hspace{0.1cm} for $(x,v) \in (x_0, \infty) \times [0,1]$.
\item[\emph{(ii)}] $\big |\partial_xB(x,v)\big | ~\lesssim~ \frac{\psi(x)}{x}v^2$  \hspace{0.1cm} for $(x,v) \in (x_0, \infty) \times [0,1]$
\item[\emph{(iii)}]   $\big |\partial_xB(x,v) \big |~\lesssim~ \frac{\psi(x)}{x} e^{\mathrm cv}$ \hspace{0.1cm} for $(x,v) \in (x_0, \infty) \times [0,\infty)$.
\end{enumerate}
\end{lemma}

\bigskip

\textbf{Proof of Lemma \ref{lem:propA}.} (i) By definition of $A$, for all $x>x_{\psi}$, $v \geq 0$
$$
\frac{|A(x,v)|}{x^{\frac{\kappa}{1-\kappa}}} \leq e^{v \frac{\kappa}{1-\kappa}} \int_x^{xe^v} \frac{|h(u)|}{u^{\frac{\kappa}{1-\kappa}}} \mathrm du \leq e^{v \frac{\kappa}{1-\kappa}} ~\bar h_{\frac{\kappa}{1-\kappa}}(x) x(e^v-1).
$$
We conclude with  $e^{v \frac{\kappa}{1-\kappa}}(e^v-1)-v \lesssim v^2$ on  $[0,1]$.

\smallskip

(ii) Just use that for all $x>x_{\psi},v \geq 0$,  $\partial_x A(x,v)=h(xe^v)e^v-h(x)$.

\smallskip

(iii) Rewrite $\partial_x A(x,v)=h(x)(e^v-1)+(h(xe^v)-h(x))e^v$. By the mean value theorem, for all $x>x_{\psi}, v>0$, there exists $c_{x,v} \in (x,xe^v)$ such that
$$
\partial_x A(x,v)=h(x)(e^v-1)+h'(c_{x,v})x(e^v-1)e^v,
$$
and then
$$
\frac{\partial_x A(x,v)}{x^{\frac{2\kappa}{1-\kappa}}}\leq \frac{h(x)}{x^{\frac{2\kappa}{1-\kappa}}}(e^v-1)+\frac{h'(c_{x,v})}{(c_{x,v})^{\frac{2\kappa}{1-\kappa}}}e^{v^{\frac{2\kappa}{1-\kappa}}}x(e^v-1)e^v,
$$
which gives the result since  $e^v-1 \lesssim v$ and $(e^v-1)e^{v(1+\frac{2\kappa}{1-\kappa})} -v \lesssim v^2$ when $v \in [0,1]$.
$\hfill \square$
\bigskip

\textbf{Proof of Lemma \ref{lem:propB}.} By the mean value theorem, for all $u>x >x_{\psi}$, 
$$
\psi(u)=\psi(x)+\psi'(x) (u-x)+\frac{1}{2} \psi''(c_{x,u})(u-x)^2
$$
for some $c_{x,u} \in (x,u)$.

\smallskip

(i) With this notation, for all $x >x_{\psi}, v>0$,
\begin{eqnarray*}
B(x,v)&=&\int_x^{x e^v}\frac{\psi(x)+\psi'(x)(u-x)+\frac{1}{2} \psi''(c_{x,u})(u-x)^2}{u} \mathrm du-\psi(x)v \\
&=& \psi'(x) x(e^v-1-v)+ \frac{1}{2} \int_x^{x e^v}\frac{\psi''(c_{x,u})(u-x)^2}{u} \mathrm du,
\end{eqnarray*}
and therefore for $(x,v) \in (x_{\psi}, \infty) \times [0,1]$
$$
\left |B(x,v)- \psi'(x)x\frac{v^2}{2} \right| \lesssim \psi'(x) xv^3+ \sup_{u \in [x,xe^v]} |\psi''(u)|x^2(e^v-1)^3.
$$
By Lemma \ref{prop:boundpsi} there exists $x_0 \geq x_{\psi}$ and $\mathrm c \in (0,\infty)$ such that
$$
\psi'(x) \lesssim \frac{\psi(x)}{x}\quad \text{and}  \quad \big| \psi''(x)\big| \lesssim \frac{\psi(x)}{x^2} \quad \text{on } (x_0,\infty) 
$$
and
$$
\psi(xe^v) \leq e^{\mathrm cv} \psi(x) \quad \text{ for } (x,v) \in (x_0,\infty) \times [0,\infty).
$$
Hence,
\begin{eqnarray*}
\left |B(x,v)- \psi'(x)x\frac{v^2}{2} \right| &\lesssim& \psi(x) v^3+ \sup_{u \in [x,xe^v]} \frac{\psi(u)}{u^2}x^2(e^v-1)^3 \\
& \underset{\text{$\psi$ is increasing}} \lesssim &\psi(x) v^3+\psi(xe^v)v^3 \\
&\lesssim & \psi(x) v^3.
\end{eqnarray*}

\smallskip

(ii)-(iii) Apply again the mean value theorem, for all $x >x_{\psi}, v> 0$,
\begin{equation*}
\partial_x B(x,v)=\frac{\psi(xe^v)}{x}-\frac{\psi(x)}{x}-\psi'(x)v = \psi'(x) (e^v-1-v)+\frac{\psi''(c_{x,xe^v})}{2}x(e^v-1)^2.
\end{equation*}
Hence for $(x,v) \in (x_0,\infty) \times [0,\infty)$,
\begin{eqnarray*}
\left |\partial_x B(x,v) \right | &\lesssim&  \frac{\psi(x)}{x} (e^v-1-v)+ \sup_{u \in [x,xe^v]} \frac{\psi(u)}{u^2} x (e^v-1)^2 \\
&\leq&   \frac{\psi(x)}{x}  \left( (e^v-1-v) + e^{\mathrm cv} (e^v-1)^2\right)
\end{eqnarray*}
where as above we use Lemma \ref{prop:boundpsi} and that $\psi$ is increasing. This gives the expected inequalities for points (ii) and (iii).
$\hfill \square$

\subsection{Refinement of the bounds for $h,h'$}
\label{sec:refin}

We now use the rough bounds of Proposition \ref{prop_majo_hprime} and the results on $A,B$ established in the previous section to get finer (and of exact order) bounds for $h$ and $h'$, using again Lemma \ref{prop:contract} and Lemma \ref{prop:derivative}. These finer bounds will allow us to quickly prove Proposition \ref{prop:equiv}, which is done in the next section. Roughly our strategy holds as follows. Using (\ref{hAB}), we would like to set up an inequality of the form 
$$
h(x) \leq \mathrm a h(x)+ \frac{b}{x}
$$
for some  $\mathrm a \in (0,1)$ and $\mathrm b>0$ and all $x$ large enough. This would imply that $h(x) \leq \frac{b}{(1-\mathrm a)x}$ as expected. Unfortunately this is not possible to do directly, mainly because we cannot bound from above $\int_x^{x e^v}h(u) \mathrm du$ by something like $h(x)x(e^v-1)$, which would have permit us to get the expected $\mathrm a \in (0,1)$ thanks to the inequality of Lemma \ref{prop:contract}. However we are not far from being able to do it, the idea being to use the function $\bar h_{\frac{\kappa}{1-\kappa}}$ instead of $h$.  We will get the $O(x^{-2})$ bound for $h'$ similarly. 

\begin{lemma}
\label{lem:hgrandO}
$h(x)=O(x^{-1})$ as $x \rightarrow \infty$.
\end{lemma}

\medskip

\textbf{Proof.} We use the expression (\ref{hAB}) of $h$ and split the integral into two. Note first, using the positivity of $A(x,v),B(x,v)$ for $x>x_+$ and all $v\geq 0$ that
\begin{eqnarray*}
0\leq \int_1^{\infty}\exp(-\psi(x)v) \left(1-\exp\left(-(A(x,v)+B(x,v))\right) \right)\pi(\mathrm dv) &\leq& \int_1^{\infty}\exp(-\psi(x)v) \pi(\mathrm dv) \\
&=& O(\exp(-ax)) \text{ as } x \rightarrow \infty
\end{eqnarray*}
for some $a>0$ by Lemma \ref{prop:derivative} (iv) and since $x=O(\psi(x))$ as  $x \rightarrow \infty$.
Then, use that $\exp(-y) \leq y$, $\forall y \geq 0$, together with Lemma \ref{lem:propA} (i) and Lemma \ref{lem:propB} (i) to get the existence of positive constants $b,c$ such that for all $x$ large enough
\begin{eqnarray*}
0 &\leq& \int_0^{1}\exp(-\psi(x)v) \left(1-\exp\left(-(A(x,v)+B(x,v))\right) \right)\pi(\mathrm dv) \\
&\leq & \int_0^{1}\exp(-\psi(x)v) \left(A(x,v)+B(x,v) \right)\pi(\mathrm dv) \\
&\leq & \bar h_{\frac{\kappa}{1-\kappa}}(x) x^{1+\frac{\kappa}{1-\kappa}} \int_0^{1} \exp(-\psi(x)v) (v+b v^2) \pi(\mathrm dv) + c \psi(x)  \int_0^{1}\exp(-\psi(x)v) v^2\pi(\mathrm dv) \\
&\leq & \kappa \bar h_{\frac{\kappa}{1-\kappa}}(x)  x^{\frac{\kappa}{1-\kappa}} + \frac{bc x^{\frac{\kappa}{1-\kappa}} }{\psi(x)} +\frac{bc}{x},
\end{eqnarray*}
the last inequality being a consequence of Lemma \ref{prop:contract} and Lemma \ref{prop:derivative} (iii).
Summing the two inequalities to get an upper bound for $h$ (using its expression (\ref{hAB})) and recalling that $\psi(x)\rightarrow \infty$ as $x \rightarrow \infty$, we therefore see that there exists some constants $\mathrm a \in (\kappa,1)$ and $\mathrm b \in (bc,\infty)$ such that for all $x$ large enough,
$$
0 \leq h_{\frac{\kappa}{1-\kappa}}(x)  \leq \mathrm a \bar h_{\frac{\kappa}{1-\kappa}}(x) + \frac{\mathrm b}{x^{1+\frac{\kappa}{1-\kappa}}}.
$$
The function in the right-hand side of this inequality being decreasing, we get that for all $x$ large enough 
$$
0 \leq \bar h_{\frac{\kappa}{1-\kappa}}(x)  \leq \sup_{y\geq x} \left(\mathrm a \bar h_{\frac{\kappa}{1-\kappa}}(y) + \frac{\mathrm b }{y^{1+\frac{\kappa}{1-\kappa}}}\right)=\mathrm a \bar h_{\frac{\kappa}{1-\kappa}}(x) + \frac{\mathrm b }{x^{1+\frac{\kappa}{1-\kappa}}}
$$
and therefore, since $ \bar h_{\frac{\kappa}{1-\kappa}}(x)$ is finite,
$$
0 \leq \bar h_{\frac{\kappa}{1-\mathrm \kappa}}(x)  \leq \frac{\mathrm b }{(1-\mathrm a)x^{1+\frac{\kappa}{1-\kappa}}}
$$
which finally leads to
$$
 0 \leq h(x) \leq \frac{\mathrm b}{(1-\mathrm a)x}
$$
for all $x$ large enough. $\hfill \square$

\bigskip

\smallskip

Since $A(x,v)=\int_x^{xe^v} h(u) \mathrm du$, Lemma \ref{lem:hgrandO} immediately yields:

\begin{cor}
\label{cor:majoAprem}
\begin{enumerate}
\item[\emph{(i)}] $|A(x,v)| \lesssim v \quad \text{ for } (x,v) \in (x_{\psi}+1,\infty) \times [0,1].$
\item[\emph{(ii)}] $|\partial_xA(x,v)| \lesssim x^{-1} \quad \text{ for } (x,v) \in (x_{\psi}+1,\infty) \times [0,\infty).$   
\end{enumerate}
\end{cor}

Then, using this result, we prove with arguments similar to those used for Lemma  \ref{lem:hgrandO} that:

\begin{lemma}
\label{lem:hprimegrandO}
$h'(x)=O(x^{-2})$ as $x \rightarrow \infty$.
\end{lemma}

\medskip

\textbf{Proof of Lemma \ref{lem:hprimegrandO}.} 
Thanks to Lemma \ref{lem:propA}, Lemma \ref{lem:propB}, Proposition \ref{prop_majo_hprime} and Lemma \ref{prop:derivative} we can use the Dominated Convergence Theorem to derive (\ref{hAB}) and get for $x$ large enough
\begin{eqnarray}
\label{h'AB}
h'(x) &=& -\psi'(x) \int_0^{\infty} \exp(-\psi(x)v) \left(1-\exp\left(-(A(x,v)+B(x,v))\right) \right)\ v \pi(\mathrm dv)  \\
&+& \int_0^{\infty}\exp(-\psi(x)v) \exp\left(-(A(x,v)+B(x,v))\right)\left( \partial_xA(x,v)+\partial_xB(x,v)\right)\pi(\mathrm dv) \nonumber.
\end{eqnarray}

\smallskip

$\bullet$ For $x$ large enough the first line in the right-hand side of (\ref{h'AB}) is negative and its absolute value is smaller than
\begin{eqnarray*}
&& \psi'(x)  \int_1^{\infty} \exp(-\psi(x)v) \left(1-\exp\left(-(A(x,v)+B(x,v))\right) \right) v \pi(\mathrm dv) \\
&+&\psi'(x) \int_0^{1} \exp(-\psi(x)v) \left(A(x,v)+B(x,v)\right) v \pi(\mathrm dv) \\
&\leq& \psi'(x)  \int_1^{\infty} \exp(-\psi(x)v) v \pi(\mathrm dv) \\
&+& a\psi'(x) \left( \int_0^1 \exp(-\psi(x)v) v^2 \pi(\mathrm dv) + \psi(x)  \int_0^1 \exp(-\psi(x)v)\right) v^3 \pi(\mathrm dv),
\end{eqnarray*}
where $a$ denotes some positive constant (independent of $x$): this a consequence of the positivity of $A(x,v),$ $B(x,v)$ for $x$ large enough and all $v$, of Corollary \ref{cor:majoAprem} (i) and Lemma \ref{lem:propB} (i). The terms \linebreak $\int_0^1 \exp(-\psi(x)v) v^2 \pi(\mathrm dv)$, $ \psi(x)  \int_0^1 \exp(-\psi(x)v)) v^3 \pi(\mathrm dv)$ and $\int_1^{\infty} \exp(-\psi(x)v) v \pi(\mathrm dv)$ are all smaller than a constant times $(x\psi(x))^{-1}$ by Lemma \ref{prop:derivative}. So finally the absolute value of the first line in the right-hand side of (\ref{h'AB}) is smaller than a constant times
$\frac{ \psi'(x)}{x \psi(x)} = O(x^{-2})$ as $x \rightarrow \infty$, the $O(x^{-2})$ being a consequence of Lemma \ref{prop:boundpsi}. 

\smallskip

$\bullet$ The absolute value of the second line in the right-hand side of (\ref{h'AB})  is smaller than
$$
 \int_0^{\infty}\exp(-\psi(x)v)\left(| \partial_x A(x,v)| + |\partial_x B(x,v)| \right)\pi(\mathrm dv).
$$
By Lemma \ref{lem:propA} (iii), Lemma \ref{lem:hgrandO} and Lemma \ref{prop:contract}, Lemma \ref{lem:propB} (ii) and Lemma \ref{prop:derivative} (iii), there exists a constant $b$ such that for all $x$ large enough
$$
 \int_0^{1}\exp(-\psi(x)v)\left( | \partial_x A(x,v)| + |\partial_x B(x,v)| \right)\pi(\mathrm dv) \leq  \bar h'_{\frac{2\kappa}{1-\kappa}} (x) x^{1+\frac{2\kappa}{1-\kappa}}\left( \frac{\kappa}{x} +\frac{b}{x\psi(x)}\right)+O(x^{-2}).
$$
On the other hand, by Corollary \ref{cor:majoAprem} (ii), Lemma \ref{lem:propB} (iii) and Lemma \ref{prop:derivative} (iv):
$$
 \int_1^{\infty}\exp(-\psi(x)v))\left(| \partial_x A(x,v)| + |\partial_x B(x,v)| \right)\pi(\mathrm dv) = O\left(\exp(-cx)\right) \quad \text{as } x \rightarrow \infty
 $$
for some $c>0$. 

\smallskip

$\bullet$ Gathering these informations (and recalling that $\psi(x) \rightarrow \infty$ as $x \rightarrow \infty$), we have shown  that there exists some constants $\mathrm a \in (\kappa,1)$ and $\mathrm b \in (0\infty)$ such that for all $x$ large enough
$$
\frac{|h'(x)|}{x^{\frac{2\kappa}{1-\kappa}}} \leq \mathrm a  \bar h'_{\frac{2\kappa}{1-\kappa}} (x) + \frac{\mathrm b}{x^{2+\frac{2\kappa}{1-\kappa}}}.
$$
The function in the right-side being decreasing, this implies that for $x$ large enough
$$
\bar h'_{\frac{2\kappa}{1-\kappa}} (x) \leq  \mathrm a  \bar h'_{\frac{2\kappa}{1-\kappa}} (x) + \frac{\mathrm b}{x^{2+\frac{2\kappa}{1-\kappa}}}
$$
and then, recalling that $ \bar h'_{\frac{2\kappa}{1-\kappa}} (x) $ is finite (by Lemma \ref{prop_majo_hprime}), that
$$
\bar h'_{\frac{2\kappa}{1-\kappa}} (x) \leq \frac{\mathrm b}{(1-\mathrm a)x^{2+\frac{2\kappa}{1-\kappa}}}
$$
which leads to the expected $h'(x)=O(x^{-2})$ as $x \rightarrow \infty$.
 $\hfill \square$
 
 \bigskip
 
This leads to a more precise control of the function $A$:
 
\begin{cor}
\label{cor:majoAdeu}
$$
\left| A(x,v) -h(x)x v \right| \lesssim v^2 \quad \text{for } (x,v) \in (x_{\psi+1},\infty) \times [0,1]. 
$$
\end{cor}

\textbf{Proof.} By the mean value theorem, for $x > x_{\psi},v>0$,
\begin{eqnarray*}
A(x,v)=\int_x^{xe^v}h(u) \mathrm du=\int_x^{xe^v}\left(h(x)+h'(c_{u,x})(u-x) \right)\mathrm du
\end{eqnarray*}
where $c_{u,x} \in (x,u)$ for all $u \in (x,xe^v)$.
Hence
$$
\left| A(x,v) -h(x)x v \right| \leq  h(x)x(e^v-1-v) +\frac{1}{2}\max_{u \in [x,xe^v]}|h'(u)| x^2 (e^v-1)^2
$$
and we conclude with Lemma \ref{lem:hgrandO} and Lemma \ref{lem:hprimegrandO}, and the fact that $(e^v-1-v)\lesssim v^2,$ $(e^v-1)^2 \lesssim v^2$ for $v \in [0,1]$.
$\hfill \square$

\subsection{Proof of  Proposition \ref{prop:equiv}}
\label{sec:proofProp}

We repeat the arguments of the proof of Lemma \ref{lem:hgrandO}  with the finer bounds obtained in the previous section which allows us to obtain the expected asymptotics. As already noticed, as $x \rightarrow \infty$,
$$
h(x)=\int_0^1 \exp(-\psi(x)v) \left(1-\exp\left(-A(x,v)-B(x,v)\right) \right)\pi(\mathrm dv)+ O(\exp(-ax))
$$ 
for some $a>0$. Then, since $y-y^2 \leq 1-\exp(-y) \leq y$ for $y\geq 0$ and $(y_1+y_2)^2 \leq 2 (y_1^2+y_2^2)$ for $(y_1,y_2) \in \mathbb R^2$, we see that for all $x$ large enough
\begin{eqnarray*}
&& \left |\int_0^1 \exp(-\psi(x)v) \left(1-\exp\left(-(A(x,v)+B(x,v))\right) \right)\pi(\mathrm dv) \right. \\
&& \qquad \left.-\int_0^{1} \exp(-\psi(x) v)  \big(A(x,v)+B(x,v)\big)\pi (\mathrm dv) \right | \\
&\leq& \mathrm 2   \int_0^{1} \exp(-\psi(x) v) \big(A^2(x,v)+B^2(x,v)\big) \pi (\mathrm dv) \\
&=& O\left( \frac{1}{x \psi( x)}\right) \quad \text{as }x \rightarrow \infty
\end{eqnarray*}
where for the last inequality we use Corollary \ref{cor:majoAprem} (i), Lemma \ref{lem:propB} (i) and Lemma \ref{prop:derivative} (iii). 

Besides, by Corollary \ref{cor:majoAdeu} and Lemma \ref{prop:derivative} (i) and (iii),
\begin{eqnarray*}
\int_0^{1} \exp(-\psi(x) v)  A(x,v) \pi (\mathrm dv)&=&h(x) x \frac{1}{x} \left(1- \frac{\psi(x)}{x\psi'(x)}\right) +O\left(\frac{1}{x \psi( x)}\right) \\
&=& h(x) -h(x)  \frac{\psi(x)}{x\psi'(x)} + O\left(\frac{1}{x \psi( x)}\right)
\end{eqnarray*}
and by Lemma \ref{lem:propB} (i) and Lemma \ref{prop:derivative} (ii) and (iii),
\begin{eqnarray*}
\int_0^{1} \exp(-\psi(x) v)  B(x,v) \pi (\mathrm dv)&=& \psi'(x) x \frac{1}{2} \left( \frac{2}{x^2 \psi'(x)} \left(1- \frac{\psi(x)}{x\psi'(x)}\right)-\frac{\psi(x)\psi''(x)}{x^2(\psi'(x))^3} \right)+O\left(\frac{1}{x \psi( x)}\right) \\
&=&   \frac{1}{x} \left(1- \frac{\psi(x)}{x\psi'(x)}\right)-\frac{\psi(x)\psi''(x)}{2x(\psi'(x))^2} +O\left(\frac{1}{x \psi( x)}\right)
\end{eqnarray*}
So finally, as $x \rightarrow \infty$,
\begin{eqnarray*}
h(x)&=& \int_0^{1} \exp(-\psi(x) v)  \big(A(x,v)+B(x,v)\big)\pi (\mathrm dv) +O  \left( \frac{1}{x \psi( x)}\right)\\
&=& h(x)-h(x)  \frac{\psi(x)}{x\psi'(x)} + \frac{1}{x} \left(1- \frac{\psi(x)}{x\psi'(x)}\right)-\frac{\psi(x)\psi''(x)}{2x(\psi'(x))^2} +O\left(\frac{1}{x \psi( x)}\right) 
\end{eqnarray*}
and therefore 
\begin{eqnarray*}
h(x) &=& \frac{x\psi'(x)}{\psi(x)} \left(\frac{1}{x} \left(1- \frac{\psi(x)}{x\psi'(x)}\right)-\frac{\psi(x)\psi''(x)}{2x(\psi'(x))^2} +O\left(\frac{1}{x \psi( x)}\right) \right)\\
&=& \frac{\psi'(x)}{\psi(x)} -\frac{1}{x} - \frac{\psi''(x)}{2 \psi'(x)} +O\left(\frac{\psi'(x)}{ (\psi( x))^2}\right). 
\end{eqnarray*}

\section{A common case: Laplace exponents with asymptotic power expansions}
\label{sec:special}

In this section we give more concrete expressions of the equivalences of Corollary \ref{cor:RV}  and discuss applications when the function $\phi$ has the following asymptotic expansion (it is always implied that $I$ is the exponential functional of a subordinator with L\'evy measure $\pi$ our Laplace exponent $\phi$): 

\begin{lemma}
\label{prop:aenuinfinite}
Assume that $\phi$ writes
$$\phi(x)=x^{\gamma} \left(1-\sum_{i=1}^p c_i{x^{-\gamma_i}}+O(x^{-1-\varepsilon})\right) \text{ as } x \rightarrow \infty$$
for some $\gamma \in [0,1)$, $1/2<\gamma_1 <\gamma_2 \ldots <\gamma_{p-1} < \gamma_p = 1$, $\varepsilon>0$ and $c_i \in \mathbb R$, $1 \leq i \leq p$ . Then, 
$$
\frac{\psi(x)}{x}=x^{\frac{\gamma}{1-\gamma}}-\sum_{i=1}^p \frac{c_i}{(1-\gamma)}x^{\frac{\gamma-\gamma_i}{1-\gamma}}+O(x^{-1-\eta})
$$
for some $\eta>0$, and consequently 
$$
\mathbb P(I>t) ~\underset{t \rightarrow \infty}\sim ~ k(t) t^{-\frac{\gamma}{1-\gamma}} \underset{t \rightarrow \infty}\propto t^{\frac{c_p-\frac{\gamma}{2}}{(1-\gamma)}}\exp\left(- (1-\gamma)t^{\frac{1}{1-\gamma}}+\sum_{i=1}^{p-1} \frac{c_i}{1-\gamma_i}  t^{\frac{1-\gamma_i}{1-\gamma}}\right).
$$
\end{lemma}

\bigskip

\textbf{Proof.}  
By definition of $\psi$ and using the assumption on $\phi$, we get for all $x>0$
$$
x=\frac{\psi(x)}{\phi(\psi(x))}=\frac{\psi(x)^{1-\gamma}}{ \left(1-\sum_{i=1}^p c_i\psi(x)^{-\gamma_i}+O(\psi(x)^{-1-\varepsilon}) \right)}.
$$
This implies that $\psi(x) \underset{x \rightarrow \infty}\sim x^{\frac{1}{1-\gamma}}$ and more precisely that
\begin{eqnarray*}
\psi(x) &=& x^{\frac{1}{1-\gamma}} \left(1-\sum_{i=1}^p c_i\psi(x)^{-\gamma_i}+O\left(\psi(x)^{-1-\varepsilon}\right)  \right)^{\frac{1}{1-\gamma}} \\
&=& x^{\frac{1}{1-\gamma}} \left(1-\sum_{i=1}^p \frac{c_i}{(1-\gamma)}\psi(x)^{-\gamma_i}+O(\psi(x)^{-2\gamma_1})+O\left(\psi(x)^{-1-\varepsilon}\right) \right).
\end{eqnarray*}
In particular $\psi(x)^{\gamma_i}=x^{\frac{\gamma_i}{1-\gamma}}\left(1+O(\psi(x)^{-\gamma_1})\right)$, $1 \leq i \leq p$. By reinjecting this back into the above asymptotic expansion of $\psi$ we get 
\begin{eqnarray*}
\psi(x)= x^{\frac{1}{1-\gamma}}-\sum_{i=1}^p \frac{c_i}{(1-\gamma)}x^{\frac{1-\gamma_i}{1-\gamma}}+O\Big(x^{\frac{1-2\gamma_1}{1-\gamma}}\Big)+O\Big(x^{-\frac{\varepsilon}{1-\gamma}}\Big).
\end{eqnarray*}
So if $2\gamma_1>1$, then $\eta:=\min(2\gamma_1-1,\varepsilon)>0$ and
$$
\frac{\psi(x)}{x}=x^{\frac{\gamma}{1-\gamma}}-\sum_{i=1}^p \frac{c_i}{(1-\gamma)}x^{\frac{\gamma-\gamma_i}{1-\gamma}}+O(x^{-1-\eta}).$$
We conclude with Corollary \ref{cor:RV}.
$\hfill \square$

\bigskip

If now some $\gamma_i$s in the asymptotic expansion of $\phi$ are smaller than $1/2$ they will also contribute to the non-$O(x^{-1-\eta})$ terms of $\psi(x)/x$ which in turn will contribute to additional terms in or in front of the exponential. This quickly gives complex formulas when several $\gamma_i$s are close to 0.

\bigskip

\textbf{Applications.} This allows us in particular to proove the examples discussed in Section \ref{sec:MR}:

\medskip

\textbf{1. Case 1 and Case 2 of Section \ref{sec:MR}.} Using the definition of $\phi$ and an integration by part, the above lemma easily leads to the results stated there.

\medskip

\textbf{2. ``$\mathbf{(a,b,c)}$-L\'evy measures".} To treat this family of examples, we could either use the two cases mentioned above or directly that:
\begin{lemma}
\label{lem:abBeta}
Let $a>0,b>-1$. Then
\begin{eqnarray*}
\frac{1}{\Gamma(b)}\int_0^1(1-u^x)u^{a-1}(1-u)^{b-1} \mathrm du &=& \frac{\Gamma(a)}{\Gamma(a+b)}- \frac{\Gamma(x+a)}{\Gamma(x+a+b)} \\
&\underset{\text{as $x \rightarrow \infty$}}{=}& \frac{\Gamma(a)}{\Gamma(a+b)}-  x^{-b}+b\left(a+\frac{b}{2}-\frac{1}{2}\right)x^{-b-1}+O(x^{-b-2})
\end{eqnarray*}
where $1/\Gamma$ denotes the extension to $\mathbb C$ by analytic continuation of the function $1/\Gamma$ initially defined on $\{z\in \mathbb C:\mathrm{Re}(z)>0\}$.
\end{lemma}

\textbf{Proof.} Fix $x\geq 1$ and $a>0$, and let $z \in \mathbb C$, $\mathrm{Re}(z)>-1$. If $\mathrm{Re}(z)>0$, we can split the integral
\begin{eqnarray*}
&& \frac{1}{\Gamma(z)}\int_0^1(1-u^x)u^{a-1}(1-u)^{z-1} \mathrm du \\
&=&\frac{1}{\Gamma(z)}\int_0^1u^{a-1}(1-u)^{z-1} \mathrm du-\frac{1}{\Gamma(z)}\int_0^1u^xu^{a-1}(1-u)^{z-1} \mathrm du \\
&=&  \frac{\Gamma(a)}{\Gamma(a+z)}- \frac{\Gamma(x+a)}{\Gamma(x+a+z)}.
\end{eqnarray*}
Then by analytic continuation, the identity extends to all $z \in \mathbb C$ such that $\mathrm{Re}(z)>-1$. Next, fix $a>0,b>-1$, apply the above identity with $z=b$ and use the asymptotic expansion as $x \rightarrow \infty$ of order 2 of the Gamma function
$$
\Gamma(x)=e^{-x} x^{x-1/2} \sqrt{2\pi} \left( 1+\frac{1}{12x}+O(x^{-2})\right)
$$
to get for all $c>0$
$$\frac{\Gamma (x+c)}{\Gamma(x)}= x^{c}\left(1-\frac{c(1-c)}{2x}+O(x^{-2}) \right),$$
which leads to the asymptotic expansion of the statement. $\hfill \square$

\bigskip

Combining this asymptotic expansion with Lemma \ref{prop:aenuinfinite} allows us to illustrate the variety of situations:

\begin{lemma}
\label{lem:detailabc}
Assume that $\pi(\mathrm dx)= c^{-1}e^{-ax}(1-\exp(-x/c))^{b-1}\mathrm dx$, with $a,c>0, b>-1$. Then:
\begin{itemize}
\item[$\bullet$] if $b \in (-1,-1/2)$, 
\begin{eqnarray*}
\mathbb P(I>t) ~ &\underset{t \rightarrow \infty}\sim ~& k(t) t^{\frac{b}{1+b}} \\
&\underset{t \rightarrow \infty} \propto & t^{\left(ba+\frac{b(b-1)}{2c} +\frac{b}{2}\right)\frac{1}{1+b}} \exp\left( -(1+b)\left(\frac{|\Gamma(b)|}{c^b} \right)^{\frac{1}{1+b}}t^{\frac{1}{1+b}}+ \frac{|\Gamma(b)| \Gamma(ac)}{\Gamma(b+ac)(1+b)} t\right)
\end{eqnarray*}
\item[$\bullet$] if $b \in (1/2,1)$, 
$$
\mathbb P(I>t) ~\underset{t \rightarrow \infty}\sim ~ k(t) \underset{t \rightarrow \infty}\propto \exp\left(-\frac{\Gamma(ac)\Gamma(b)}{\Gamma(ac+b)}t +\left( \frac{\Gamma(ac+b)}{\Gamma(ac)\Gamma(b) c}\right)^b \frac{\Gamma(b)}{1-b} t^{1-b}\right)
$$
\item[$\bullet$] if $b=1$,
$$
\mathbb P(I>t) ~\underset{t \rightarrow \infty}\sim ~ k(t) \underset{t \rightarrow \infty}\propto t^{a} \exp(-t/ac).
$$
\item[$\bullet$] if $b>1$,
$$
\mathbb P(I>t) ~\underset{t \rightarrow \infty}\sim ~ k(t) \underset{t \rightarrow \infty}\propto  \exp\left(-\frac{\Gamma(ac)\Gamma(b)}{\Gamma(ac+b)}t\right).
$$
\end{itemize}
For $b \in [-1/2,1/2]\backslash\{0\}$ there will be more and more terms in or in front of the exponential as it approaches 0.
\end{lemma}

\textbf{Proof.} For $b \neq 0$, just note that 
\begin{eqnarray}
\label{eqq:phi}
\nonumber
\phi(x)&=&\int_0^1 (1-v^{xc})v^{ac-1}(1-v)^{b-1} \mathrm dv \\
&=& \Gamma(b) \left( \frac{\Gamma(ac)}{\Gamma(ac+b)}- \frac{\Gamma(xc+ac)}{\Gamma(xc+ac+b)} \right) \\
\nonumber
&\underset{\text{as $x \rightarrow \infty$}}{=}& \Gamma(b) \left(\frac{\Gamma(ac)}{\Gamma(ac+b)}-  (cx)^{-b}+b\left(ac+\frac{b}{2}-\frac{1}{2}\right)(cx)^{-b-1}+O(x^{-b-2}) \right)
\end{eqnarray}
by Lemma  \ref{lem:abBeta}. Then use Lemma \ref{prop:aenuinfinite} to conclude, noticing that the dominant term in the asymptotic expansion of $\phi$ depends on the sign of $b$. The details of the calculations are easy and left to the reader. $\hfill \square$

\bigskip

\textbf{Remark.} In this model, there are two particular cases where the distribution of $I$ is explicit:
\begin{enumerate}
\item[$\bullet$]
When $a=1$ and $b=-c \in (-1,0)$, the random variable $I$ is proportional to a generalized Mittag-Leffler random variable with parameters $(c,c)$. We recall that a random variable $M$ is said to have a generalized Mittag-Leffler distribution with parameters $(\alpha,\theta)$, with $\alpha \in (0,1)$, $\theta>-\alpha$ if for all suitable test functions $f$,
$$
\mathbb E[f(M)]=C \mathbb E[\sigma_{\alpha}^{-\theta}f(\sigma_{\alpha}^{-\alpha})]
$$
where $\sigma_{\alpha}$ is a stable random variable with Laplace transform $\mathbb E[\exp^{-q\sigma_{\alpha}}]=\exp(-q^{\alpha})$, $q\geq 0$, and $C$ the appropriate normalizing constant. It has positive moments of all orders given by
$$
\mathbb E[M^r]=\frac{\Gamma(\theta)\Gamma(\theta/\alpha+r)}{\Gamma(\theta/\alpha)\Gamma(\theta+r \alpha)}, \quad \text{for all }r \geq 0.
$$ 
See e.g. Section 0.5 of Pitman's book \cite{PitmanStFl} for details and for a series expansion of the density of such a random variable when $\theta=0$, from which one easily deduce the general $\theta$-case.
\item[$\bullet$] For any $a>0$, if $c=(a+1)^{-1}$ and $b=-1+(a+1)^{-1}$, then $I$ is proportional to $\mathbf e(1)^{1/(a+1)}$, where $\mathbf e(1)$ denotes an exponential random variable with parameter 1.  
\end{enumerate}
To see these two identities in distribution, it suffices to compute the entire moments of $I$ (which characterize its distribution). On the one hand, it is known from Carmona, Petit, and Yor \cite[Prop]{CPY97} that
$$
\mathbb E[I^n]=\frac{n!}{\prod_{i=1}^n \phi(i)}.
$$
On the other hand we see by (\ref{eqq:phi}) that when $ac+b=0$ (and since $1/\Gamma(0)=0$), the function $\phi$ is proportional to $\Gamma((x+a)c)/\Gamma(xc)$. Together with the Carmona, Petit, Yor's expression of the moments of $I$ we see that when moreover $a=1$ (hence $b=-c$), the $n-$th moment of $I$ is proportional to
$$
\frac{n! \Gamma(c)}{\Gamma((n+1)c)}
$$
with indeed coincides with the $n-$th moment generalized Mittag-Leffler random variable with parameters $(c,c)$. Whereas when $c=(a+1)^{-1}$ (hence $b=-a(a+1)^{-1}$), whatever $a>0$, the $n-$th moment of $I$ is proportional to (using that $i \Gamma(\frac{i}{a+1}) = (a+1)\Gamma(\frac{i}{a+1}+1)$)
$$
\frac{n! }{\prod_{i=1}^n \frac{\Gamma(\frac{i+a}{a+1})}{\Gamma(\frac{i}{a+1})}}=(a+1)^n \prod_{i=1}^n \frac{\Gamma(\frac{i+a+1}{a+1})}{\Gamma(\frac{i+a}{a+1})}=(a+1)^n \Gamma\left(\frac{n}{a+1}+1\right)
$$
which is the $n$-th moment of the random variable $(a+1)\mathbf e(1)^{1/(a+1)}$.

\section{Application to self-similar Markov processes}
\label{sec:ssMp}

This short section is devoted to the proof of Corollary \ref{cor:ssMp}. We use the notations introduced in the neighborhood of this corollary, in particular we recall that $X^{(1)}=\exp(-\xi_{\rho(\cdot)})$ is a non-increasing non-negative $1/\alpha-$self-similar Markov process that starts from 1 and reaches 0 at time
$$
I_X=\int_0^{\infty}\exp(-\alpha \xi_r) \mathrm dr,
$$
with $\xi$ a drift-free and unkilled subordinator with Laplace exponent $\phi$, and $\psi$ the inverse of $x\mapsto x/\phi(x)$. Observe that $\psi(\alpha \cdot)/\alpha$ is the inverse of $x\mapsto x/\phi(\alpha x)$, with $\phi(\alpha \cdot)$ the Laplace exponent of $\alpha\xi$. Corollary \ref{cor:ssMp} is, in fact, nearly a direct consequence of the combination of our Theorem \ref{thm:main} with Theorem 1.2 of \cite{HR12} which states that under (\ref{hyp:main}) (see also \cite[Theorem 3.1]{H10} under the more restrictive assumption that $\phi$ is regularly varying):
$$
\frac{\psi(\alpha t)}{\alpha t}(X^{(1)}(t))^{\alpha}  \text{ conditioned on }  X^{(1)}(t)>0 \text{ converges in distribution as }t\rightarrow \infty \text{ to } R,
$$
where $R$ is a positive random variable whose distribution is characterized by the fact that when $R$ is taken independently of $I_X$, $RI_X$ has an exponential distribution with parameter 1  (the existence and uniqueness of such a random variable is proved in \cite{BY01}; it is characterized by its entire moments: $\mathbb E[R^n]=\prod_{i=1}^n \phi(\alpha i)$). To get  Corollary \ref{cor:ssMp}, we need however to improve the above conditional convergence by showing that the moments of any positive order also converge. 

\bigskip

For this we use arguments already used here or there. 
To lighten notations, we write $X$ instead of $X^{(1)}$ in the few lines that follow. The key observation for the proofs of \cite[Theorem 1.2]{HR12} and \cite[Theorem 3.1]{H10} is that $(I_X-t)^+$ is distributed as $X^{\alpha}(t)\tilde I_X$ where $\tilde I_X$ is independent of $X(t)$ and distributed as $I_X$, for all $t\geq 0$ (this is immediate, using the definition of $I_X$ and the independence and stationarity of the increments of $\xi$). An integration by parts then gives for $a>0$
$$
\mathbb E\left[X^a(t)\right] =\frac{\frac{a}{\alpha}\int_0^{\infty} u^{\frac{a}{\alpha}-1} \mathbb P(I_X>u+t) \mathrm du}{\mathbb E\Big[I_X^{\frac{a}{\alpha}}\Big]}
$$
and then
\begin{equation}
\label{eq:dom}
\left(\frac{\psi(\alpha t)}{\alpha t}\right)^{\frac{a}{\alpha}} \frac{\mathbb E\left[X^a(t)\right]}{\mathbb P(I_X>t)} =\frac{\frac{a}{\alpha}}{\mathbb E\Big[I_X^{\frac{a}{\alpha}}\Big]} \int_0^{\infty} u^{\frac{a}{\alpha}-1} \frac{\mathbb P\left(I_X>u \frac{\alpha t}{\psi(\alpha t)}+t\right)}{\mathbb P(I_X>t)} \mathrm du.
\end{equation}
Besides, from Proposition 3.1 of \cite{HR12}, under (\ref{hyp:main}),
$$
\frac{\mathbb P\left(I_X>u\frac{\alpha t}{\psi(\alpha t)}+t\right)}{\mathbb P(I_X>t)} \underset{t\rightarrow \infty} \longrightarrow \exp(-u), \quad \forall u \geq 0
$$
and we claim that for all $u>0$ and all $t$ large enough (independent of $u$) 
$$
\frac{\mathbb P(I_X>u\frac{\alpha t}{\psi(\alpha t)}+t)}{\mathbb P(I_X>t)} \leq \exp(-u).
$$
Indeed, if we set $f(x):=-\ln (\mathbb P(I_X>x))$ we know from Proposition \ref{lem:hpositive}  that $f'(x) \geq \frac{\psi(\alpha x)}{\alpha x} $ for all $x$ large enough. Using that $x \mapsto \frac{\psi(x)}{x}$ is non-decreasing, this implies that for $t$ large enough and all $u>0$
$$
f\left(t+u\frac{\alpha t}{\psi(\alpha t)}\right) - f(t)=\int_t^{t+u\frac{\alpha t}{\psi(\alpha t)}} f'(x) \mathrm dx \geq \frac{\psi(\alpha t)}{\alpha t} \left(t+u\frac{\alpha t}{\psi(\alpha t)}-t\right)=u
$$
and therefore 
$$
\frac{\mathbb P(I_X>t+u\frac{\alpha t}{\psi(\alpha t)})}{\mathbb P(I_X>t)} =\exp\left(-f\left(t+u\frac{\alpha t}{\psi(\alpha t)}\right)+f(t)\right) \leq \exp(-u).
$$
We can therefore apply the dominated convergence theorem to get from (\ref{eq:dom}) that
$$
\left(\frac{\psi(\alpha t)}{\alpha t}\right)^{\frac{a}{\alpha}} \frac{\mathbb E\left[X^a(t)\right]}{\mathbb P(I_X>t)}  \underset{t \rightarrow \infty}\longrightarrow \frac{\frac{a}{\alpha}}{\mathbb E\Big[I_X^{\frac{a}{\alpha}}\Big]} \int_0^{\infty} u^{\frac{a}{\alpha}-1} \exp(-u)\mathrm du=\mathbb E[R^{\frac{a}{\alpha}}].
$$
Combined with Theorem \ref{thm:main}, this gives Corollary \ref{cor:ssMp}.

\bigskip

\bibliographystyle{siam}
\bibliography{frag}

\end{document}